%% file: agt-5-44.tex
\documentclass{gtart_h}
\input agtout

\lognumber{44}
\volumenumber{5}
\volumeyear{2005}
\papernumber{44}
\pagenumbers{1075}{1109}
\received{26 October 2004} 
\accepted{28 June 2005}
\published{31 August 2005}

\usepackage{graphicx,psfrag}
\usepackage{amsmath,amssymb,mathrsfs}
\usepackage{enumerate}        
\usepackage[all]{xy}    
\UseComputerModernTips        

\def\psfraga <#1,#2> #3#4{%
\psfrag {#3}{\smash{\rlap{\kern #1 \raise #2\hbox{#4}}}}}

\def\figref#1{\hyperlink{#1anchor}{Figure~\ref*{#1}}}
\def\anchor#1{\noindent\hypertarget{#1anchor}{\smash{$\phantom{99}$}}}

\newcommand{\im}{\operatorname{im}}
\newcommand{\dom}{\operatorname{dom}}

\newcommand {\bng}{B_n\Gamma}
\newcommand {\cng}{\mathcal{C}^n\Gamma}
\newcommand {\ucng}{U\mathcal{C}^n\Gamma}
\newcommand {\dng}{\mathcal{D}^n\Gamma}
\newcommand {\udng}{U\mathcal{D}^n\Gamma}

\newtheorem{theorem}{Theorem}[section]
\newtheorem{lemma}[theorem]{Lemma}
\newtheorem{proposition}[theorem]{Proposition}
\newtheorem{corollary}[theorem]{Corollary}

\theoremstyle{definition}

\newtheorem{example}[theorem]{Example}

\begin{document}

\title{Discrete Morse theory and graph braid groups}
\authors{Daniel Farley\\Lucas Sabalka}
\address{Department of Mathematics, University of Illinois at 
               Urbana-Champaign\\Champaign, IL  61820, USA}
 \asciiemail{farley@math.uiuc.edu, sabalka@math.uiuc.edu}
 \gtemail{\mailto{farley@math.uiuc.edu}, \mailto{sabalka@math.uiuc.edu}}

\gturl{\url{www.math.uiuc.edu/~farley}{\rm\qua and\qua}\url{www.math.uiuc.edu/~sabalka}}
\asciiurl{www.math.uiuc.edu/ farley, www.math.uiuc.edu/ sabalka}

\begin{abstract}
If $\Gamma$ is any finite graph, then the \emph{unlabelled configuration
space of $n$ points on $\Gamma$}, denoted $\ucng$, is the space of
$n$-element subsets of $\Gamma$.  The \emph{braid group of $\Gamma$ on $n$
strands} is the fundamental group of $\ucng$.

We apply a discrete version of Morse theory to these $\ucng$, for any $n$
and any $\Gamma$, and provide a clear description of the critical cells in
every case.  As a result, we can calculate a presentation for the braid
group of any tree, for any number of strands.  We also give a simple proof
of a theorem due to Ghrist: the space $\ucng$ strong deformation retracts
onto a CW complex of dimension at most $k$, where $k$ is the number of
vertices in $\Gamma$ of degree at least $3$ (and $k$ is thus independent
of $n$).
\end{abstract}
\asciiabstract{%
If Gamma is any finite graph, then the unlabelled configuration space
of n points on Gamma, denoted UC^n(Gamma), is the space of n-element
subsets of Gamma.  The braid group of Gamma on n strands is the
fundamental group of UC^n(Gamma).  We apply a discrete version of
Morse theory to these UC^n(Gamma), for any n and any Gamma, and
provide a clear description of the critical cells in every case.  As a
result, we can calculate a presentation for the braid group of any
tree, for any number of strands.  We also give a simple proof of a
theorem due to Ghrist: the space UC^n(Gamma) strong deformation
retracts onto a CW complex of dimension at most k, where k is the
number of vertices in Gamma of degree at least 3 (and k is thus
independent of n).}

\primaryclass{20F65, 20F36}
\secondaryclass{57M15, 57Q05, 55R80}
\keywords{Graph braid groups, configuration spaces, discrete Morse theory}

\maketitle


\section{Introduction}

Let $\Gamma$ be a finite graph, and fix a natural number $n$.  The
\emph{labelled configuration space $\cng$} is the $n$-fold Cartesian
product of $\Gamma$, with the set $\Delta = \{ (x_1, \ldots, x_n) \mid x_i
= x_j \hbox{ for some } i \neq j \}$ removed.  The \emph{unlabelled
configuration space $\ucng$} is the quotient of $\cng$ by the action of
the symmetric group $S_{n}$, where the action permutes the factors.  The
fundamental group of $\ucng$ (respectively, $\cng$) is the \emph{braid
group} (respectively, the \emph{pure braid group}) \emph{of $\Gamma$ on
$n$ strands}, denoted $B_{n}\Gamma$ (respectively, $PB_{n}\Gamma$).

Configuration spaces of graphs occur naturally in certain motion planning
problems.  An element of $\cng$ or $\ucng$ may be regarded as the
positions of $n$ robots on a given track (graph) within a factory.  A path
in $\cng$ represents a collision-free movement of the robots. Ghrist and
Koditschek \cite{GK} used a repelling vector field to define efficient and
safe movements for the fundamentally important case of two robots on a
Y-shaped graph.  Farber \cite{Fa} has recently shown that $TC(\cng)=2k+1$,
where ``$TC$" denotes topological complexity, $\Gamma$ is any tree, $k$ is
the number of vertices of $\Gamma$ having degree at least $3$, and $n \geq
2k$.

In \cite{G}, Ghrist showed that $\cng$ strong deformation retracts onto a
complex $X$ of dimension at most $k$, where $k$ is the number of vertices
having degree at least three in $\Gamma$ (and thus is independent of $n$).  
If $\Gamma$ is a \emph{radial tree}, i.e., if $\Gamma$ is a tree having
only one vertex of degree more than 2, then $\cng$ strong deformation
retracts on a graph.  By computing the Euler characteristic of this graph,
Ghrist computes the rank of the pure braid group on $\Gamma$ as a free
group. (See also \cite{Ga}, where the Euler characteristic of any
configuration space of a simplicial complex is computed.)

Ghrist also conjectures in \cite{G} that every braid group $\bng$ is a
\emph{right-angled Artin group} -- i.e.\ a group which has a presentation in
which all defining relations are commutators of the generators.  Abrams
\cite{A1} revised the conjecture to apply only to planar graphs. These
conjectures were one of our original reasons to investigate graph braid
groups.  However, Abrams \cite{A3} has shown that some graph braid groups,
including a planar graph and a tree braid group, are not Artin
right-angled.  There are positive results: Crisp and Wiest \cite{CW} have
recently shown that every graph braid group embeds in a right-angled Artin
group.

We use Forman's discrete Morse theory \cite{Fo} (see also \cite{B}) in
order to simplify any space $\ucng$ within its homotopy type.  We are able
to give a very clear description of the critical cells of a discretized
version of $\ucng$ with respect to a certain discrete gradient vector
field $W$ (see Section 2 for definitions).  By a theorem of Forman's
(\cite{Fo}, page 107), a CW complex $X$ endowed with a discrete gradient
vector field $W$ is homotopy equivalent to a complex with $m_{p}$ cells of
dimension $p$, where $m_{p}$ is the number of cells of dimension $p$ that
are critical with respect to $W$.  (A similar theorem was proved by Brown
earlier; see page 140 in \cite{B}.)

Our classification of critical cells leads to a simple proof of Ghrist's
theorem that the braid group of any radial tree is free, and we compute
the rank of this free group as a function of the degree of the central
vertex and the number of strands.  This computation features an unusual
application of the formula for describing the number of ways to distribute
indistinguishable balls into distinguishable boxes.  We also prove a
somewhat strengthened version of Ghrist's strong deformation retraction
theorem by a simple application of Morse-theoretic methods.  Our dimension
bounds resemble those of \cite{Sw}, where the homological dimension of the
braid groups $\bng$ was estimated.  The strong deformation retract in our
theorem has an explicit description in terms of so-called collapsible,
redundant, and critical cells (see Section 2 for definitions).
        
The main theorem of the paper describes how to compute a presentation of
$\bng$, where $n$ is an arbitrary natural number and $\Gamma$ is any tree.  
We deduce our theorem with the aid of a particular discrete gradient
vector field $W$ and the ``redundant 1-cells lemma", which greatly
simplifies calculations.  The generators in our presentation correspond to
critical $1$-cells of $\ucng$, and the relators correspond to the critical
$2$-cells.  The relators are all commutators, so that the rank of the
abelianization of $\bng$ is equal to the number $m_{1}$ of critical
$1$-cells.  It follows easily that any CW complex homotopy equivalent to
$\ucng$ must have at least $m_{1}$ $1$-cells.  It seems reasonable to
guess that, more generally, any CW complex homotopy equivalent to $\ucng$
must have at least $m_{p}$ $p$-cells, where $m_{p}$ is the number of
$p$-cells that are critical with respect to $W$, although we don't prove
this.  The presentations that we obtain for $\ucng$ do not look like
presentations of right-angled Artin groups, but we don't know how to show
that tree braid groups are not right-angled Artin groups by our methods.
  
It is possible to calculate group presentations for any graph braid group
using the techniques of this paper, but we offer no general theorem here,
since the resulting presentations are somewhat less than optimal -- e.g.,
our calculations in some cases yield presentations of the free group which
contain non-trivial relators.  We hope to improve these presentations in 
the near future.
               
Our presentation is mostly self-contained.  Section 2 contains a short
exposition of Forman's Morse theory, which is sufficient for all of the
applications in later sections.  The central idea here is that of a
discrete gradient vector field, which induces a classification of the
cells into one of three mutually exclusive types:  collapsible, redundant,
and critical.  Section 3 describes how to define a discrete gradient
vector field on the configuration space of any graph, and gives a
description of the collapsible, redundant, and critical cells with respect
to the given discrete gradient vector field.  Section 4 contains a
calculation of the braid group of a radial tree, and a refined version of
Ghrist's strong deformation retraction theorem.  Section 5 contains the
main theorem, about presentations of tree braid groups, with examples.

We would like to thank Ilya Kapovich and Robert Ghrist for participating
in numerous discussions relating to this work.  We thank Kim Whittlesey,
Aaron Abrams, James Slack and the referee for reading an earlier version
of the manuscript, and suggesting corrections.  The first author thanks
Tadeusz Januszkiewicz for telling him of the references \cite{Ga} and
\cite{Sw}.


\section{Preliminary material}



\subsection{Discrete Morse theory}


The fundamental idea of discrete Morse theory is that of a    
\emph{collapse}.  We take the definition from \cite{C}, page 14.  
If $(X,Y)$ is a finite CW pair then \emph{$X$ collapses to $Y$ by an 
elementary collapse} -- denoted $X \searrow\!\!\!\!^{e}\,\, Y$ -- if and only if

\begin{enumerate}
\item $X = Y \cup e^{n-1} \cup e^{n}$ where $e^{n}$ and $e^{n-1}$ are    
open cells of dimension $n$ and $n-1$, respectively, which are not in    
$Y$, and
\item there exists a ball pair $(Q^{n},Q^{n-1})\approx (I^{n},I^{n-1})$    
and a map $\phi: Q^{n} \rightarrow X$ such that

\begin{enumerate}
\item $\phi$ is a characteristic map for $e^{n}$
\item $\phi \mid Q^{n-1}$ is a characteristic map for $e^{n-1}$
\item $\phi(P^{n-1}) \subseteq Y^{n-1}$, where $P^{n-1} = cl(\partial 
Q^{n} - Q^{n-1})$.
\end{enumerate}
\end{enumerate}
We say that \emph{$X$ collapses to $Y$}, and write $X \searrow Y$, if  
$Y$ may be obtained from $X$    
by a sequence of elementary collapses.      

Let $X$ be a finite CW complex.  Let $K$ denote the set of open cells of
$X$ with $K_p$ the set of open $p$-dimensional cells of $X$.  For an open
cell $\sigma \in K$, we write $\sigma^{(p)}$ to indicate that $\sigma$ is
of dimension $p$.  We write $\sigma < \tau$ if $\sigma \neq \tau$ and
$\sigma \subseteq \overline{\tau}$ (where $\overline{\tau}$ denotes the
closure of $\tau$).  We write $\sigma \leq \tau$ if $\sigma = \tau$ or
$\sigma < \tau$.

We will need to work with a special type of CW complex.  From now on,
every CW complex $X$ we consider will have the following property: if
$\sigma^{(p)}$ and $\tau^{(p+1)}$ are open cells of $X$ and $\sigma^{(p)}
\cap \overline{\tau^{(p+1)}} \neq \emptyset$, then $\sigma^{(p)} <
\tau^{(p+1)}$.  The importance of this assumption will become apparent
in Proposition \ref{prop:filtration}.

Suppose that $\sigma^{(p)}$ is a face of $\tau^{(p+1)}$ ($\sigma <
\tau$).  Let $B$ be a closed ball of dimension $p+1$, and let $h: B \to
X$ be a characteristic map for $\tau$. The cell $\sigma$ is a
\emph{regular face} of $\tau$ if $h: h^{-1}\sigma \to \sigma$ is a
homeomorphism, and $\overline{h^{-1}(\sigma)}$ is a closed $p$-ball.

A \emph{discrete vector field} $W$ on $X$ is a sequence of partial
functions $W_i: K_i \to K_{i+1}$ such that:

\begin{enumerate}[(i)]
\item Each $W_i$ is injective
\item If $W_i(\sigma) = \tau$, then $\sigma$ is a regular face of $\tau$.
\item $\im (W_i) \cap \dom (W_{i+1}) = \emptyset$.
\end{enumerate}

This definition of a discrete vector field differs very slightly from
Forman's (\cite{Fo}, page 130-131).  Note that a \emph{partial function}
from a set $A$ to a set $B$ is simply a function defined on a subset of
$A$.

Let $W$ be a discrete vector field on $X$.  A \emph{$W$-path of dimension
$p$} is a sequence of $p$-cells $\sigma_0, \sigma_1, \dots, \sigma_r$ such
that if $W(\sigma_i)$ is undefined, then $\sigma_{i+1} = \sigma_i$, and
otherwise $\sigma_{i+1} \neq \sigma_i$ and $\sigma_{i+1} < W(\sigma_i)$.
The $W$-path is \emph{closed} if $\sigma_r = \sigma_0$, and
\emph{non-stationary} if $\sigma_1 \neq \sigma_0$. A discrete vector field
$W$ is a \emph{discrete gradient vector field} if $W$ has no
non-stationary closed paths.

Given any discrete gradient vector field $W$ on $X$, there is an
associated classification of cells in $X$ into 3 types: redundant,
collapsible, and critical (this terminology is partially borrowed from
\cite{Fo} as well as from Ken Brown, \cite{B}).  A cell $\sigma \in K$ is
\emph{redundant} if $\sigma \in \dom W$, \emph{collapsible} if $\sigma \in
\im W$, and \emph{critical} otherwise. The \emph{rank} of a cell $c$ with
respect to a discrete gradient vector field $W$ is the length of the
longest $W$-path $c = c_{1},\ldots,c_{r}$ having the property that $c_{i}
\neq c_{j}$ if $i \neq j$.  Critical and collapsible cells all have rank
1, and redundant cells are of rank at least 2.  If $c' < W(c)$ and $c \neq
c'$, then clearly $rank(c')<rank(c)$.

In \cite{Fo} (page 131), Forman shows that discrete gradient vector fields
are precisely the discrete vector fields which arise from discrete Morse
functions (in a manner he describes).  We will work directly with discrete
gradient vector fields, and never with an explicit discrete Morse
function.


\subsection{Monoid presentations and rewrite systems} 
\label{sec:monoids&rewriting}


An \emph{alphabet} is simply a set $\Sigma$.  The \emph{free monoid on
$\Sigma$}, denoted $\Sigma^{\ast}$, is the collection of all positive
words in the generators $\Sigma$, together with the empty word, endowed
with the operation of concatenation.

A \emph{monoid presentation}, denoted $\langle \Sigma \mid \mathcal{R}
\rangle$, consists of an alphabet $\Sigma$ together with a collection
$\mathcal{R}$ of ordered pairs of elements in $\Sigma^{\ast}$.  An
element of $\mathcal{R}$ should be regarded as an equality between words
in $\Sigma^{\ast}$, but, in what follows, the order will matter.

A \emph{rewrite system $\Gamma$} is an oriented graph.  The vertices of
$\Gamma$ are called \emph{objects} and the positive edges are called
\emph{moves}.  If $v_{1}$ is the initial vertex of some positive edge in
$\Gamma$ and $v_{2}$ is the terminal vertex, then write $v_{1}
\rightarrow_{\Gamma} v_{2}$, or $v_{1} \rightarrow v_{2}$ if the name of
the specific rewrite system is clear from the context.  An object is
called \emph{reduced} if it is not the initial vertex of any positive
edge (move).  The reflexive, transitive closure of $\rightarrow$ is
denoted $\dot{\rightarrow}$.

A rewrite system is called \emph{terminating} if every sequence $a_{1}
\rightarrow a_{2} \rightarrow a_{3} \rightarrow \ldots$ is finite.  A
rewrite system is called \emph{confluent} if, whenever $a
\dot{\rightarrow} b$ and $a \dot{\rightarrow} c$, there is an object $d$
such that $b \dot{\rightarrow} d$ and $c \dot{\rightarrow} d$. A rewrite
system is \emph{locally confluent} if when $a \rightarrow b$ and $a
\rightarrow c$, then there is $d$ such that $b \dot{\rightarrow} d$ and
$c \dot{\rightarrow} d$.

\begin{lemma}\label{Nwmn} {\rm\cite{N}}\qua    
Every terminating locally confluent rewrite system is confluent.
\end{lemma}

A rewrite system is called \emph{complete} if it is both terminating and
confluent.  For a complete rewrite system, it is not difficult to argue
that every equivalence class of the equivalence relation generated by
$\rightarrow$ has a unique reduced object.

Every monoid presentation $\langle \Sigma \mid \mathcal{R} \rangle$ has a
natural rewrite system, called a \emph{string rewriting system},
associated to it.  The set of objects of this string rewriting system is
the free monoid $\Sigma^{\ast}$.  There is a move from $w_{1} \in
\Sigma^{\ast}$ to $w_{2} \in \Sigma^{\ast}$ if $w_{1} = u r_{1} v$ and
$w_{2} = u r_{2} v$ in $\Sigma^{\ast}$, where $u, v \in \Sigma^{\ast}$,
$(r_{1}, r_{2}) \in \mathcal{R}$.


\subsection{Discrete Morse theory and group presentations}\label{sec:GP}


Assume in this section that $X$ is a finite connected CW complex with a
discrete gradient vector field $W$. Let $X_{n}'$ be the subcomplex of $X$
consisting of the $n$-skeleton, but with the redundant $n$-cells removed.  
Let $X_{n}''$ consist of the $n$-skeleton of $X$, but with the redundant
and critical $n$-cells removed.  The following proposition was essentially
proved by Brown (\cite{B}, page 140) in the case of semi-simplicial
complexes.

\begin{proposition} \label{prop:filtration}

Consider the following filtration of $X$:
  $$ \emptyset = X_{0}'' \subseteq X_{0}' \subseteq \ldots \subseteq    
  X_{n}'' \subseteq X_{n}' \subseteq X_{n+1}'' \subseteq \ldots.$$     

\begin{enumerate}
\item For any $n$, $X_{n}'$ is obtained from $X_{n}''$ by attaching    
$m_{n}$ $n$-cells to $X_{n}''$ along their boundaries, where $m_{n}$ is    
the number of critical $n$-cells of the discrete gradient vector field    
$W$.
\item For any $n$, $X_{n+1}'' \searrow X_{n}'$.
\end{enumerate}
\end{proposition}

\begin{proof}(1)\qua This is obvious.

(2)\qua Let $X_{n,k}$ be the subcomplex of $X$ consisting of the entire
$(n-1)$-skeleton, together with all $n$-cells of rank at most $k$, and all
$(n+1)$-cells that are the images of such $n$-cells under the function
$W$.  Thus, for example, $X_{n}' = X_{n,1}$.

We claim that $X_{n,i+1} \searrow X_{n,i}$, for
$i \in \mathbb{N}$.

Let $c$ be an open $n$-cell of rank exactly $i+1$.  Since $W$ is
injective, $W(c)$ cannot be the image under $W$ of a cell of rank less
than or equal to $i$, and so it lies outside of $X_{n,i}$.  If $c \cap
X_{n,i}$ is nonempty, it can only be because there is some open
collapsible $(n+1)$-cell $c'$ such that $c \cap \overline{c'}$ is nonempty
and $c' = W(c'')$, for some open $n$-cell $c''$ of rank less than or equal
to $i$.  Given our standing assumption about the CW complex $X$ (from
Subsection \ref{sec:monoids&rewriting}), we thus know that $c < c'$.  Now
if $c = c_1, c_2, c_3, \ldots c_{i+1}$ is a $W$-path without any repeated
cells (which exists because $rank(c)=i+1$) we have that $c'', c_1, c_2,
\ldots, c_{i+1}$ is also a $W$-path without repetitions, since clearly
$c'' \neq c_1$ and there are no non-stationary closed $W$-paths.  This
implies that the rank of $c''$ is at least $i+2$, a contradiction.  It
follows that the first part of the definition of a collapse is satisfied
for the pair $(X_{n,i+1}, X_{n,i+1} - (c \cup W(c)))$.

Since $c$ is a regular face of $W(c)$, there is a characteristic map
$\phi: B^{n+1} \rightarrow X_{n,i+1}$ for $W(c)$ such that $\phi:
\phi^{-1}(c) \rightarrow c$ is an homeomorphism and
$\overline{\phi^{-1}(c)}$ is a closed ball.  It follows easily that the
second part of the definition of a collapse is satisfied for the same
pair.

Repeating this argument for every $(n+1)$-cell of $X_{n,i+1}$, we
eventually conclude that $X_{n,i+1} \searrow X_{n,i}$, since the
individual collapses are compatible.

It follows that $X_{n,k} \searrow X_{n}'$ for any $k \in \mathbb{N}$.  
For $k$ sufficiently large, $X_{n,k}$ consists of the entire $n$-skeleton,
plus the collapsible $n+1$-cells.  That is, $X_{n+1}'' = X_{n,k} \searrow
X_{n}'$.
\end{proof}

We collect a number of corollaries in the following proposition:

\begin{proposition} \label{prop:corollaries}$\phantom{99}$

\begin{enumerate}
\item The inclusion of $X_{n}'$ into $X$ induces an isomorphism from
$\pi_{n-1}(X_{n}')$ to $\pi_{n-1}(X)$.

\item If $X$ has no critical cells of dimension greater than $k$, then
$X \searrow X_{k}'$.

\item {\rm(\cite{Fo}, page 107)}\qua $X$ is homotopy equivalent to a CW complex with    
$m_{n}$ cells of dimension $n$, where $m_{n}$ is the number of critical    
$n$-cells in $X$.

\item The subcomplex of $X$ generated by the collapsible and critical
edges is connected.

\item The subcomplex of $X$ generated by the collapsible edges and the    
$0$-skeleton of $X$ is a maximal forest.

\item If there is only one critical $0$-cell, then the graph consisting    
of (the closures of) the collapsible edges is a maximal tree in $X$.
\end{enumerate}
\end{proposition}

\begin{proof}$\phantom{99}$

(1)\qua Note that the map $\pi_{n-1}(X_{n}') \rightarrow \pi_{n-1}(X)$
factors as
  $$ \pi_{n-1}(X_{n}') \rightarrow \pi_{n-1}(X_{n+1}'') \rightarrow    
  \pi_{n-1}(X).$$
The first map is an isomorphism, because $X_{n+1}'' \searrow X_{n}'$.     
The second map is also an isomorphism, since $X$ is obtained from
$X_{n+1}''$ by attaching cells of dimension greater than or equal to
$n+1$, which have no effect on $\pi_{n-1}$.

(2)\qua We have the sequence     
  $$ X_{k}' \subseteq X_{k+1}'' = X_{k+1}' \subseteq X_{k+2}'' = X_{k+2}'    
  \subseteq \ldots.$$
Each complex in this sequence collapses onto the one before it (sometimes
trivially, where equality holds).  Since the sequence terminates at $X$,
(2) is proved.

(3)\qua This follows easily from the previous proposition.

(4)\qua In fact, the subcomplex in question is $X_{1}'$, and
$\pi_{0}(X_{1}') \rightarrow \pi_{0}(X)$ is a bijection by (1).  Since
$X$ is connected, so is $X_{1}'$.

(5)\qua Since $X_{1}'' \searrow X_{0}'$, each component of $X_{1}''$ is
contractible, and so $X_{1}''$ is a forest.  It is true by definition that
$X_{1}''$ contains the whole $0$-skeleton, so $X_{1}''$ is also maximal.

(6)\qua Since $X_{1}'' \searrow X_{0}'$, and $X_{0}'$ is a singleton set,
$X_{1}''$ is connected.  By (5), we know that $X_{1}''$ is also a maximal
forest, so it must in fact be a maximal tree.
\end{proof}

Choose a maximal tree $T$ of $X$ consisting of all of the collapsible
edges in $X$, and additional critical edges, as needed.

Define a monoid presentation $\mathcal{MP}_{W,T}$ as follows: Generators
are oriented edges in $X$, both positive and negative, so that there are
two oriented edges for each geometric edge in $X$.  If $e$ denotes a
particular oriented edge, let $\overline{e}$ denote the edge with the
opposite orientation.  If $w$ denotes a sequence of oriented edges
$e_1\dots e_m$, let $\overline{w}$ denote sequence of oriented edges
$\overline{e_m} \dots \overline{e_1}$.

A \emph{boundary word} of a $2$-cell $c$ is simply one of the possible
relations determined by an attaching map for $c$ (cf.\ \cite{St}, page 139);  
if $w_{1}$ and $w_{2}$ are two boundary words for a cell $c$, then $w_{1}$
can be obtained from $w_{2}$ by the operations of inverting and taking
cyclic shifts.

There are several types of relations.

\begin{enumerate}
\item For a given oriented edge $e$ in $T$, introduce the relations    
$(e,1)$ and $(\overline{e},1)$.
\item For any oriented edge $e$, introduce relations $(e\overline{e},1)$    
and $(\overline{e}e,1)$.
\item For a collapsible $2$-cell $c$, consider the (unique) geometric    
$1$-cell $e$ such that $e =W^{-1}(c)$.  Suppose that a boundary word of    
$c$ is $ew$.  In this case, the word $w$ contains no occurrence of $e$ or    
$\overline{e}$, since the geometric edge corresponding to $e$ is a    
regular face of $c$.  Introduce the relations $(e,\overline{w})$ and    
$(\overline{e},w)$.  

\end{enumerate}

\begin{proposition} \label{cmplt} 
The rewrite system associated to the monoid presentation    
$\mathcal{MP}_{W,T}$ is complete.
\end{proposition}      

\begin{proof}

If $w \rightarrow w_{1}$ and $w \rightarrow w_{2}$ correspond to disjoint
applications of relations in $\mathcal{MP}_{W,T}$, that is, if $w =
rs_{1}tu_{1}v$, $w_{1}=rs_{2}tu_{1}v$, and $w_{2}=rs_{1}tu_{2}v$, where
$r$, $s_{1}$, $s_{2}$, $t$, $u_{1}$, $u_{2}$, and $v$ are words in the
free monoid on oriented edges in $X$, and $(s_{1},s_{2})$,
$(u_{1},u_{2})$ are relations in $\mathcal{MP}_{W,T}$, then $w_{1}
\rightarrow w_{3}$ and $w_{2} \rightarrow w_{3}$, where $w_{3} =
rs_{2}tu_{2}v$.

Thus we need only consider the cases in which the moves $w \rightarrow
w_{1}$ and $w \rightarrow w_{2}$ are not disjoint.  Checking definitions,
we see that it is not possible for a move of the first type to overlap
with a move of the third type, since the left side of a relation of type
$1$ is a word of length one involving only an edge in $T$, and the left
side of a relation of type $3$ is another word of length one, but
involving a redundant edge, and redundant edges lie outside of $T$.     
Similarly easy arguments show that the moves $w \rightarrow w_{1}$ and $w
\rightarrow w_{2}$ can only overlap if at least one of these moves is of 
the second type.

Suppose, without loss of generality, that $w = te\overline{e}v$ and $w_{1}
= tv$.  If, say, $w_{2} = t\overline{e}v$, and thus the move $w
\rightarrow w_{2}$ involves the application of a relation of type $1$,
then it must be that the edge corresponding to $e$ is in $T$, so that
$w_{2} \rightarrow w_{1}$.  If $w \rightarrow w_{2}$ involves the
application of a relation of type $3$, say $w_{2} = tewv$, then $w_{2}
\rightarrow t\overline{w}wv \dot{\rightarrow} tv = w_{1}$.

The case in which both moves $w \rightarrow w_{1}$ and $w \rightarrow
w_{2}$ involve relations of the second type is left as an easy exercise.
\end{proof}

In view of Proposition \ref{cmplt} (and the remarks after Lemma 
\ref{Nwmn}), there is a unique reduced word in each equivalence class 
modulo the presentation $\mathcal{MP}_{W,T}$.  If $w$ is any word in the 
generators of $\mathcal{MP}_{W,T}$, let $r(w)$ denote the unique reduced 
word that is equivalent to $w$.            

\begin{theorem} \label{thm:Morsepresentation}
Let $X$ be a finite connected CW complex with a discrete gradient vector    
field $W$.  Then:
  $$\pi_{1}(X) \cong \langle \Sigma \mid \mathcal{R} \rangle,$$
where $\Sigma$ is the set of positive critical $1$-cells that aren't
contained in $T$, and $\mathcal{R} = \{r(w) | w \hbox{ is the
boundary word of a critical 2-cell}\}$.
\end{theorem}

\begin{proof}
According to Proposition \ref{prop:corollaries}(1), $\pi_{1}\left( X_{2}'
\right) \cong \pi_{1}(X)$. The usual edge-path presentation of the
fundamental group (see \cite{St}, page 139) of $X_{2}'$ says that
  $$\pi_{1}\left( X_{2}' \right) \cong \langle \widehat{\Sigma}    
  \mid \widehat{\mathcal{R}} \rangle,$$
where $\widehat{\Sigma}$ is the set of positively oriented $1$-cells, and
the relations $\widehat{\mathcal{R}}$ are of three types:  first, there
are the boundary words of critical $2$-cells; second, there are the
boundary words of collapsible $2$-cells; third, there are words of length
one corresponding to edges in $T$.

Any collapsible $2$-cell is necessarily $W(c)$ for some redundant $1$-cell
$c$.  Choose $c$ to have the largest rank of all $1$-cells.  There is a
boundary word of $W(c)$ with the form $cw$, where $w$ is a word involving
only occurrences of $1$-cells having smaller rank than $c$.  In the
presentation of $\pi_{1}\left(X_{2}' \right)$, replace every occurrence of
$c$ or $\overline{c}$ other than the occurrence in the boundary word of
$W(c)$ with $\overline{w}$ or $w$, respectively. In the new presentation,
neither $c$ nor $\overline{c}$ occurs in any relation except the boundary
word of $W(c)$, where $c$ or $\overline{c}$ appears, but not both.  In
fact, the operation of replacing $c$ (respectively, $\overline{c}$) with
$\overline{w}$ (respectively, $w$) can change only the relations of the
first type, by the assumption about the rank of $c$ and because $c$ is not
collapsible.  Thus we can remove the generator $c$ from the new
presentation for $\pi_{1}\left(X_{2}' \right)$ along with the boundary
word for $W(c)$ to obtain another presentation of the same group. Notice
that the effect on the boundary words of the critical $2$-cells has been
to perform a sequence of reductions of type $3$.

One continues in the same way, inductively removing redundant $1$-cells
of the largest remaining rank, until all of the redundant $1$-cells have
been removed.

Next, remove the $1$-cells occurring in $T$ from the list of generators,
along with the corresponding relations (which are words of length one),
and remove all occurrences of the $1$-cells of $T$ from the remaining
relations. This procedure results in a presentation of the same group.
The final result after freely reducing is the presentation $\langle
\Sigma \mid \mathcal{R} \rangle$, since every alteration to the 
boundary words of the critical $2$-cells has been a move in the rewrite 
system corresponding to the monoid presentation $\mathcal{MP}_{W,T}$, the
remaining words are reduced modulo $\mathcal{MP}_{W,T}$, and the
remaining generators are the positively-oriented $1$-cells lying outside
of $T$.
\end{proof}

In case there is just one critical $0$-cell, the discrete gradient vector
field $W$ completely determines the maximal tree $T$, and we denote the
presentation $\langle \Sigma \mid \mathcal{R} \rangle$ from the previous
theorem $\mathcal{P}_{W}$, where the oriented CW complex is understood.  
The presentation $\mathcal{P}_{W}$ depends only on the choice of the
boundary words for the critical $2$-cells, since the string
rewriting system associated to the monoid presentation
$\mathcal{MP}_{W,T}$ is complete and by the previous theorem.


\section{Discrete gradient vector fields and graph braid groups}


Let $\Gamma$ be a graph, and fix a natural number $n$. The \emph{labelled
configuration space} of $\Gamma$ on $n$ points is the space
  $$\left(\prod^{n} \Gamma \right) - \Delta,$$
where $\Delta$ is the set of all points $(x_1, \dots, x_n) \in \prod^n
\Gamma$ such that $x_i = x_j$ for some $i \neq j$. The \emph{unlabelled
configuration space} of $\Gamma$ on $n$ points is the quotient of the
labelled configuration space by the action of the symmetric group $S_n$,
where the action permutes the factors. The \emph{braid group} of $\Gamma$
on $n$ strands, denoted $B_n\Gamma$, is the fundamental group of the
unlabelled configuration space of $\Gamma$ on $n$ strands. The \emph{pure
braid group}, denoted $PB_n\Gamma$, is the fundamental group of the
labelled configuration space.

The set of vertices of $\Gamma$ will be denoted by $V(\Gamma)$, and the
degree of a vertex $v \in V(\Gamma)$ is denoted $d(v)$.  If a vertex $v$
is such that $d(v) \geq 3$, $v$ is called \emph{essential}.

Let $\Delta'$ denote the union of those open cells of $\prod^n \Gamma$
whose closures intersect the diagonal $\Delta$. Let $\dng$ denote the
space $\prod^n \Gamma - \Delta'$. Note that $\dng$ inherits a CW complex
structure from the Cartesian product, and that a cell in $\dng$ has the
form $c_1 \times \dots \times c_n$ where each $c_i$ is either a vertex or
the interior of an edge whose closure is disjoint from the closure of
$c_j$ for $i \neq j$.  We also let $\udng$ denote the quotient of $\dng$
by the action of the symmetric group $S_n$ by permuting the coordinates.
Thus, an open cell in $\udng$ may be written $\{c_1, \dots, c_n\}$ where
each $c_i$ is either a vertex or the interior of an edge whose closure is
disjoint from the closure of $c_j$ for $i \neq j$.  The set notation is
used to indicate that order does not matter.

Under most circumstances, the labelled (respectively, unlabelled)     
configuration space of $\Gamma$ is homotopy equivalent to $\dng$
(respectively, $\udng$).  Specifically:

\begin{theorem} {\rm\cite{A1}}\qua
For any $n>1$ and any graph $\Gamma$ with at least $n$ vertices, the  
labelled (unlabelled) configuration space of $n$ points on $\Gamma$ strong 
deformation  retracts onto $\dng$ ($\udng$) if
\begin{enumerate}
\item each path between distinct vertices of degree not equal to $2$ passes  
through at least $n-1$ edges; and
\item each homotopically nontrivial path from a vertex to itself passes 
through at least $n+1$ edges.
\end{enumerate}
\end{theorem}
      
A graph $\Gamma$ satisfying the conditions of this theorem for a given
$n$ is called \emph{sufficiently subdivided} for this $n$.  It is clear
that every graph is homeomorphic to a sufficiently subdivided graph, no
matter what $n$ may be.  Throughout the rest of the paper, we work
exclusively with the space $\udng$ where $\Gamma$ is sufficiently
subdivided for $n$.  Also from now on, ``edge''  and ``cell'' will refer    
to closed objects.

We define a discrete gradient vector field $W$ using a maximal tree and
specific order on the vertices of $\Gamma$.  Choose a maximal tree $T$ in
$\Gamma$. Edges outside of $T$ are called \emph{deleted edges}.  Pick a
vertex $\ast$ of valence $1$ in $T$ to be the root of $T$.  Define an
operation $\wedge$ on vertices of $T$ to take two vertices $v_1$ and $v_2$
and yield the vertex of $T$ which is the endpoint of the geodesic segment
$[\ast,v_1] \cap [\ast,v_2]$ other than $\ast$.

For a given vertex $v \in \Gamma$, call the edges adjacent to $v$ the
\emph{directions from $v$}.  For each $v$, fix a total ordering of the
directions from $v$ by labelling each direction with a number between $0$
and $d(v) - 1$, assigning the number $0$ to the direction leading from
$v$ back to $\ast$.  (The single direction from $\ast$ is given the
number $1$.)  Also, define a function $g: V(\Gamma) \times V(\Gamma)
\rightarrow \mathbb{Z}$ such that $g(v_1,v_2)$ is the label of the
direction from $v_1$ that lies on the unique geodesic connecting $v_1$ to
$v_2$ within $T$, or $g(v_1, v_2)=0$ if $v_1 = v_2$.

Now order the vertices of $T$ as follows.  Let $v_1$ and $v_2$ be two    
vertices, and define the vertex $v_3 := v_1 \wedge v_2$.  Then $v_1 \leq    
v_2$ if and only if $v_3 = v_1$, or $v_3 \neq v_1$ and $g(v_3,v_1) <    
g(v_3,v_2)$.

For a given edge $e$, let $\iota(e)$ and $\tau(e)$ denote the endpoints of
$e$, where $\iota(e) \geq \tau(e)$ in the ordering on vertices.  For a
vertex $v \neq \ast$ we let $e(v)$ denote the (unique) edge of $T$
such that $\iota(e(v)) = v$ -- i.e.\ the edge of $T$ incident with $v$ and 
closest to $\ast$.  The following lemma is easy to prove:

\begin{lemma}[Order Lemma]
The ordering $\leq$ on vertices of $\Gamma$ is a linear order, with the    
following additional properties:
\begin{enumerate}[\rm(i)]
\item if $v_2 \in [\ast,v_1]$ then $v_2 \leq v_1$
\item if $v \in V(\Gamma)$ and $e$ is an edge of $T$ such that $e(v)    
  \cap e = \tau(e)$ and $v < \iota(e)$, then $\tau(e)$ is an    
  essential vertex and $0 < g(\tau(e),v) < g(\tau(e),\iota(e))$ (and thus
  $\tau(e) < v < \iota(e)$).\qed
\end{enumerate}
\end{lemma}

\begin{example}\label{exam:favtree} Consider the labelled tree $\Gamma$ in
\figref{fig:favtree}, which is sufficiently subdivided for $n=4$.  
This tree is especially interesting since it is the smallest tree for
which $B_n \Gamma$ ($n \geq 4$) appears not to be a right-angled Artin
group.  See Example \ref{exam:favtree4} for a presentation of $B_4\Gamma$,
and Example \ref{exam:Htree} for a discussion of the right-angled Artin
property.
   
\begin{figure}[ht!]\anchor{fig:favtree}
\begin{center}
\includegraphics{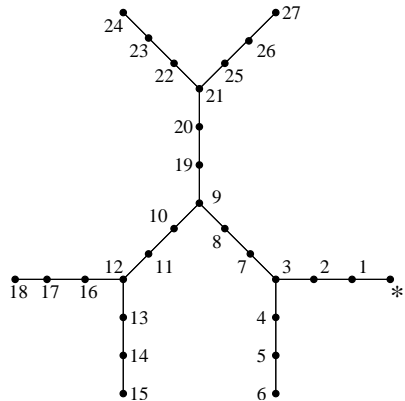}
\caption{A sufficiently subdivided tree\label{fig:favtree}}
\end{center}
\end{figure}

Let us write $v_i$ for the vertex labelled $i$.  The numbering induces an
obvious linear order on the vertices: $v_i < v_j$ if $i<j$.  This order is
a special instance of the one mentioned above, for a particular choice of
$g$ which we now describe.  It is enough to describe how to number the
directions from any given vertex $v$.  The direction from $v$ to $\ast$ is
numbered $0$, as is required, and the other directions are numbered $1, 2,
\ldots, d(v)-1$ consecutively in the clockwise order.  It is meaningful to
speak of a clockwise ordering because we have specified an embedding in
the plane by the picture.  The numbering of directions depends only on the
location of $\ast$ and the choice of the embedding.

Thus, for example, $g(v_3, v_4)=1$, $g(v_3, v_7)=2$, $g(v_9, v_{13})=1$,
and so forth.

A heuristic way to describe the numbering (or order) of the vertices is
as follows.  Begin with an embedded tree, having a specified basepoint
$\ast$ of degree one.  Now walk along the tree, following the leftmost
branch at any given intersection, and consecutively number the vertices
in the order in which they are first encountered.  (When you reach a
vertex of degree one, turn around.)  Note that this numbering depends
only on the choice of $\ast$ and the embedding of the tree.

Any ordering of the kind mentioned in the lemma may be realized by an
embedding and a choice of $\ast$, although we will not need to make use
of this fact.
\end{example}


\subsection{The function  $W$ is a discrete vector field}


Let $c = \{c_1, \dots,$ $c_{n-1},v\}$ be a cell of $\udng$ containing the
vertex $v$. If $e(v) \cap c_i = \emptyset$ for $i = 1, \dots, n-1$, then
define the cell $\{c_1, \dots, c_{n-1},e(v)\} \subset \udng$ to be the
\emph{elementary reduction} of $c$ from $v$, and we say $v$ is
\emph{unblocked} in $c$.  Otherwise, there exists some $c_i \in c$ with
$e(v) \cap c_i \neq \emptyset$, and we say $v$ is \emph{blocked by
$c_i$ in $c$}. If $v$ is the smallest unblocked vertex of $c$ in the    
sense of the order on vertices, then the reduction from $v$ is    
\emph{principal}.
        
Define a function $W$ on $\udng$ inductively.  If $c$ is a $0$-cell, let
$W_0(c)$ be the principal reduction if it exists.  For $i > 0$, let
$W_i(c)$ of an $i$-cell $c$ be its principal reduction if it exists
\emph{and} $c \not\in \im W_{i-1}$, and undefined otherwise.

Let $c$ be a cell, and let $e \in c$ be an edge of $\Gamma$.  The edge
$e$ is said to be \emph{order-respecting in $c$} if $e \subseteq T$ and,
for every vertex $v \in c$, $e(v) \cap e = \tau(e)$ implies that $v >
\iota(e)$ in the order on vertices.

For example, consider the tree in Example \ref{exam:favtree} 
(\figref{fig:favtree}).  For any vertex $v_i$, let $e_i := e(v_i)$.  For
instance, $e_{19}$ connects $v_{19}$ to $v_9$.  In the $1$-cell $\{
v_{10}, e_{19}, v_{12}, v_{16} \}$, $e_{19}$ is not order-respecting,
since $e_{10} \cap e_{19} = v_9 = \tau(e_{19})$.  On the other hand,
$e_{19}$ is order-respecting in $\{ e_{19}, v_{20}, v_{21}, v_{22} \}$.
The edge $e_{10}$ is order-respecting in any cell of $C^{4}\Gamma$.

An edge $e$ that is order-respecting in $c$ is the \emph{minimal
order-respecting edge in $c$} if $\iota(e)$ is minimal in the order on
vertices among the initial vertices of the order-respecting edges in $c$.

\begin{lemma} [Order-respecting edges lemma]$\phantom{99}$

\begin{enumerate}[\rm(i)]
\item If $\{ c_{1}, \ldots, c_{n-1}, e \}$ is obtained from      
$\{ c_{1}, \ldots, c_{n-1}, \iota(e) \}$ by principal      
reduction, then $e$ is an order-respecting edge.

\item Let $c = \{ c_{1}, \ldots, c_{n-1},e \}$ where $e$ is an edge 
contained in $T$.  
Then an edge $e'$ in $\{ c_{1}, \ldots, c_{n-1}, \iota(e) \}$ is 
order-respecting if and only if it is order-respecting in $c$.
\end{enumerate}
\end{lemma}

\begin{proof}$\phantom{99}$

(i)\qua Suppose that $e$ is not an order-respecting edge.  Clearly $e
\subseteq T$, so it must be that there is $v \in \{ c_{1}, \ldots,$
$c_{n-1}, e \}$ such that $e(v) \cap e = \tau(e)$ and $v < \iota(e)$.     
The elementary reduction from $v$ is thus well-defined in $\{ c_{1}, 
\ldots,
c_{n-1},$ $\iota(e) \}$, so the principal elementary reduction of $\{
c_{1}, \ldots,$ $c_{n-1}, \iota(e) \}$ cannot be the elementary reduction
from $\iota(e)$, since $v < \iota(e)$. We have a contradiction.

(ii)\qua We can assume that $e'$ is contained in $T$, for otherwise $e'$ 
fails to be order-respecting in every cell, and there is nothing to 
prove.

$(\Rightarrow)$\qua If $e'$ is not order-respecting in $\{c_{1},
\ldots, c_{n-1}, e \}$, then there is $v \in \{ c_{1}, \ldots,$ $ c_{n-1}, e
\}$ such that $e(v) \cap e' = \tau(e')$ and $v < \iota(e')$. Then $v \in
\{c_{1}, \ldots, c_{n-1} \} \subseteq \{ c_{1}, \ldots, c_{n-1}, \iota(e)
\}$ and thus $e'$ is not order-respecting in $\{ c_{1}, \ldots, c_{n-1},
\iota(e)\}$.

$(\Leftarrow)$\qua Suppose without loss of generality that $e' =
c_{n-1}$.  If $e'$ is not order-respecting in $c' = \{ c_{1}, \ldots,
c_{n-2}, e', \iota(e) \}$, then there is $v \in c'$ such that $e(v) \cap
e' = \tau(e')$ and $v < \iota(e')$.  In this case $v \neq \iota(e)$ since,
otherwise, $e(v)=e$ so $e \cap e' \neq \emptyset$, a contradiction.  Thus
$v \in \{ c_{1}, \ldots, c_{n-2} \} \subseteq \{ c_{1}, \ldots, c_{n-2},
e', e \}$, and so $e'$ is not order-respecting in $\{ c_{1}, \ldots,
c_{n-1}, e \}$.
\end{proof}

\begin{lemma}[Classification lemma]$\phantom{99}$
\begin{enumerate}[\rm(i)]
\item If a cell $c$ contains no order-respecting edge, then it is critical
if every vertex of $c$ is blocked, and redundant otherwise.

\item Suppose $c$ contains an order-respecting edge, and let $e$ denote
the minimal order-respecting edge in $c$.  If there is an unblocked vertex
$v \in c$ such that $v < \iota(e)$, then $c$ is redundant.  If there is no
such vertex, then $c$ is collapsible.
\end{enumerate}
\end{lemma}

\begin{proof}$\phantom{99}$

(i)\qua The previous lemma implies that $c$ cannot be in the image of $W$.     
If every vertex in $c$ is blocked, then the principal elementary reduction
of $c$ is undefined, and thus $c$ is critical. Otherwise, the principal
elementary reduction of $c$ is defined, and $c$ is redundant.

(ii)\qua Suppose first that there is an unblocked vertex $v \in c$ such that
$v < \iota(e)$.  We claim that $c = \{ c_{1}, \ldots, c_{n} \}$ is not in
the image of $W$.  If it is, then there is some $e'$ (without loss of
generality, $e' = c_{n}$) such that $\{ c_{1}, \ldots, c_{n-1}, e' \}$ is
the principal elementary reduction of $\{ c_{1}, \ldots,c_{n-1},$ 
$\iota(e') \}$.  It follows that $e'$ is an order-respecting edge in $\{
c_{1}, \ldots, c_{n-1}, e' \}$, so that $\iota(e') \geq \iota(e) > v$.     
Suppose, again without loss of generality, that $v = c_{n-1}$.  Now, since
elementary reduction from $v$ is defined for $\{ c_{1}, \ldots, c_{n-2},
v, \iota(e') \}$, it follows that $\{ c_{1}, \ldots, c_{n-1}, e' \}$ is
not the principal elementary reduction of $\{ c_{1}, \ldots, c_{n-1},
\iota(e') \}$.  This is a contradiction. Since $c$ is not in the image of
$W$ and it has unblocked vertices, it must be redundant.

Now suppose that there is no unblocked vertex $v \in c$ satisfying $v <
\iota(e)$.  Suppose, without loss of generality, that $c = \{c_1, \dots,
c_{n-1},e\}$. We claim that $\{ c_{1}, \ldots, c_{n-1}, e \}$ is the
principal reduction from $\{ c_{1}, \ldots, c_{n-1},$
$ \iota(e) \} )$.  Let
$v$ be the smallest unblocked vertex of $\{ c_{1}, \ldots, c_{n-1},
\iota(e) \}$. Clearly $ v \leq \iota(e)$, since $\iota(e)$ is unblocked.     
If $v < \iota(e)$, then $v$ is unblocked in $\{ c_{1}, \ldots, c_{n-1},
\iota(e) \}$ but blocked in $\{ c_{1}, \ldots, c_{n-1}, e \}$.  This can
only be because $e(v) \cap e = \tau(e)$. Since $e$ is order-respecting, we
have $v > \iota(e)$, which is a contradiction. Thus the principal
elementary reduction of $\{ c_{1}, \ldots, c_{n-1}, \iota(e) \}$ is $\{
c_{1}, \ldots, c_{n-1}, e\}$.

It remains to be shown that $\{ c_{1}, \ldots, c_{n-1}, \iota(e) \}$ is
not collapsible.  If it is, then there is some edge $e' \in \{ c_{1},
\ldots, c_{n-1}, \iota(e) \}$ (without loss of generality, $e' = c_{n-1}$)
such that $\{ c_{1}, \ldots, c_{n-2}, e', \iota(e) \}$ is the principal
elementary reduction of $\{ c_{1}, \ldots, c_{n-2}, \iota(e'), \iota(e)
\}$. Since $\iota(e')$ and $\iota(e)$ are both unblocked in $\{ c_{1},
\ldots, c_{n-2}, \iota(e'), \iota(e) \}$, $\iota(e') < \iota(e)$.  We
claim that $e'$ is an order-respecting edge of $\{ c_{1}, \ldots, c_{n-2},
e', e \}$. Certainly $e' \subseteq T$.  If $e'$ is not an order-respecting
edge of $\{ c_{1}, \ldots, c_{n-2}, e', e \}$, then there is some vertex
$v \in \{c_{1}, \ldots, c_{n-2}, e', e \}$ (without loss of generality, $v
= c_{n-2}$)  such that $e(v) \cap e' = \tau(e')$ and $v < \iota(e') <
\iota(e)$.  Thus $v$ is unblocked in $\{ c_{1}, \ldots, c_{n-3}, v,
\iota(e'), \iota(e) \}$ and $v < \iota(e')$, so that the elementary
reduction of $\{ c_{1}, \ldots, c_{n-3}, v, \iota(e'), \iota(e) \}$ from
$\iota(e')$ is not principal, a contradiction.  This proves that 
$e'$ is an order-respecting edge of $\{ c_{1}, \ldots , c_{n-2}, e', e \}$.

We now reach the contradiction that $e$ is not the minimal order-respecting
edge of $c$.  This completes the proof.
\end{proof}

\begin{theorem}[Classification Theorem]$\phantom{99}$
\label{thm:classification}
\begin{enumerate}       
\item A cell is critical if and only if it contains no order-respecting    
edges and all of its vertices are blocked.
\item A cell is redundant if and only if
\begin{enumerate}[\rm(a)]
\item it contains no order-respecting edges and at least one of its    
vertices is unblocked OR
\item it contains an order-respecting edge (and thus a minimal    
order-respecting edge $e$) and there is some unblocked vertex $v$ such    
that $v < \iota(e)$.
\end{enumerate}
\item  A cell is collapsible if and only if it contains an    
order-respecting edge (and thus a minimal order-respecting edge $e$) and,    
for any $v < \iota(e)$, $v$ is blocked.
\end{enumerate}
\end{theorem}

\begin{proof}
This follows logically from the previous lemma.
\end{proof}

\begin{theorem}
The function $W$ is one-to-one.
\end{theorem}

\begin{proof}
Suppose that $c = \{c_{1}, \ldots, c_{n}\}$ is collapsible.  Thus there is
a minimal order-respecting edge $e$. Since $c$ is in the image of $W$,
there must exist some edge $e'$ (without loss of generality, assume 
$c_n = e'$) such that
  $$W(\{ c_{1}, \ldots, c_{n-1}, \iota(e')\}) = \{c_{1}, \ldots, c_{n-1}, 
  e'\}.$$
A previous lemma implies that $e'$ is order-respecting in $c$.  We 
claim that $e' = e$.

If not, then $\iota(e') > \iota(e)$, $\{c_{1}, \ldots c_{n-1},
\iota(e')\}$ is redundant, and $e$ is order-respecting in $\{ c_{1},
\ldots, c_{n-1}, \iota(e')\}$. By the previous theorem, there is an
unblocked vertex $v \in \{ c_{1}, \ldots c_{n-1}, \iota(e') \}$ such that
$v < \iota(e) < \iota(e')$.  Now since $v \neq \iota(e')$, $v \in c$, in
which it must be blocked, since $c$ is collapsible.  It follows that $e(v)
\cap e' = \tau(e')$ and $v < \iota(e')$.  Since $e' \subseteq T$, we have
that $e'$ is not order-respecting in $c$, a contradiction.

It now follows that $W$ is one-to-one, since we can solve for the    
preimage of any collapsible cell.
\end{proof}      


\subsection{Proof that there are no non-stationary closed $W$-paths}


For each vertex $v$ of $\Gamma$, define a function $f_{v}$ from the cells
of $\udng$ to $\mathbb{Z}$, setting $f_{v}(c)$ equal to the number of    
$c_{i} \in c$ such that $c_{i}$ is a subset of the geodesic in $T$
connecting $\ast$ to $v$.

Each function $f_{v}$ has the following properties:

\begin{enumerate}
\item For any redundant cell $c$, $f_{v}(c) = f_{v}(W(c))$.

\item If a cell $c'$ is obtained from a cell $c$ by replacing $e \subseteq 
T$ with $\iota(e)$, then $f_{v}(c') = f_{v}(c)$.

\item If a cell $c'$ is obtained from a cell $c$ by replacing $e \subseteq 
T$ with $\tau(e)$, then $f_{v}(c') = f_{v}(c)$ unless $v \wedge \iota(e) =    
\tau(e)$, in which case $f_{v}(c') = f_{v}(c) + 1$.

\item If a cell $c'$ is obtained from a cell $c$ by replacing $e
\not\subseteq T$ with $\tau(e)$, then $f_{v}(c') = f_{v}(c)$ unless $v
\wedge \tau(e) = \tau(e)$, in which case $f_{v}(c') = f_{v}(c) + 1$.

\item If a cell $c'$ is obtained from a cell $c$ by replacing $e
\not\subseteq T$ with $\iota(e)$, then $f_{v}(c') = f_{v}(c)$ unless $v
\wedge \iota(e) = \iota(e)$, in which case $f_{v}(c') = f_{v}(c) + 1$.
\end{enumerate}

\begin{theorem}
$W$ has no non-stationary closed paths.
\end{theorem}

\begin{proof}
Suppose that $\sigma_{0}, \ldots, \sigma_{r}$ is a minimal non-stationary
closed path, so that no repetitions occur among the subsequence
$\sigma_{0}, \ldots, \sigma_{r-1}$, $\sigma_0 = \sigma_r$, and $r>1$.  
Since, for any vertex $v$ in $\Gamma$, $f_{v}(\sigma_{0}) \leq
f_{v}(\sigma_{1}) \leq \ldots \leq f_{v}(\sigma_{r}) = f_{v}(\sigma_{0})$,
equality must hold throughout.  By considering different choices for $v$,
it is clear that $\sigma_{i+1}$ may not be obtained from $W(\sigma_i)$
using rules (3), (4), or (5).  Thus, $\sigma_{i+1}$ is obtained from
$W(\sigma_{i})$ by replacing some edge $e' \in W(\sigma_{i})$ with
$\iota(e')$ (and never with $\tau(e')$), where $e'$ is necessarily
contained in $T$.  Note also that each of the cells $\sigma_{0},
\sigma_{1}, \ldots, \sigma_{r} = \sigma_{0}$ must be redundant.

We claim that if $\sigma_{i+1}$ is obtained from $W(\sigma_{i})$ by
replacing some $e \in W(\sigma_{i})$ with $\iota(e)$, then $\iota(e_{i+1})
< \iota(e_{i})$, where $e_{i}$ and $e_{i+1}$ are the minimal
order-respecting edges in $\sigma_{i}$ and $\sigma_{i+1}$, respectively.  
(If $e_{i}$ or $e_{i+1}$ doesn't exist, then $\iota(e_{i})$ or
$\iota(e_{i+1})$, respectively, is $\infty$.)

Consider first the case in which $\sigma_{i}$ has no order-respecting
edges. Since $W(\sigma_{i})$ is collapsible, it has an order-respecting
edge $e$, and $W(\sigma_{i})$ must be obtained from $\sigma_{i}$ by
replacing $\iota(e)$ with $e$ (see, for instance, the proof that $W$ is
injective).  Since $\sigma_{i+1} \neq \sigma_i$ and $\sigma_{i+1}$ must be
obtained from $W(\sigma_{i})$ by replacing an edge $e'$ from
$W(\sigma_{i})$ with $\iota(e')$, then $e' \neq e$.  Thus $e$ is the only
order-respecting edge in $\sigma_{i+1}$, and so $e_{i+1} = e$ and
$\iota(e_{i+1}) = \iota(e)<\infty = \iota(e_i)$.  This establishes the
claim if $\sigma_{i}$ has no order-respecting edges.

Now suppose that $\sigma_{i}$ has a minimal order-respecting edge $e_{i}$.
Since $\sigma_{i}$ is redundant, the minimal unblocked vertex $v$ of
$\sigma_{i}$ satisfies $v < \iota(e_{i})$.  The cell $W(\sigma_{i})$ is
obtained from $\sigma_{i}$ by replacing $v$ with $e(v)$.  Since
$\sigma_{i+1} \neq \sigma_{i}$, $\sigma_{i+1}$ must be obtained from
$W(\sigma_{i})$ by replacing some edge $e' \neq e(v)$ with $\iota(e')$.  
This implies that $e(v)$ is an order-respecting edge in $\sigma_{i+1}$
since $e(v)$ is clearly order-respecting in $W(\sigma_{i})$, and thus
$\iota(e_{i+1}) \leq v < \iota(e_{i})$.  This proves the claim.

We now reach a contradiction, because $\iota(e_{0}) > \iota(e_{1}) >
\ldots > \iota(e_{r})$, but $\sigma_{0} = \sigma_{r}$.
\end{proof}      


\subsection{Visualizing critical cells of $W$}


Using Theorem \ref{thm:classification}, it is not difficult to visualize    
the critical cells of $\udng$, for any $n$ and any $\Gamma$.  If $c$ is a
critical cell in $\udng$, then every vertex $v$ in $c$ is blocked, and
every edge $e$ in $c$ has the property that either: (i) $e$ is a deleted
edge, or (ii) $\tau(e)$ is an essential vertex, and there is some vertex
$v$ of $c$ that is adjacent to $\tau(e)$ satisfying $\tau(e) < v <
\iota(e)$.

Consider $\Gamma$ the tree given in Example \ref{exam:favtree}.  A    
critical $1$-cell in $\mathcal{C}^{4}\Gamma$ is depicted in
\figref{fig:critcell}:

\begin{figure} [ht!]\anchor{fig:critcell}
\begin{center}
\includegraphics{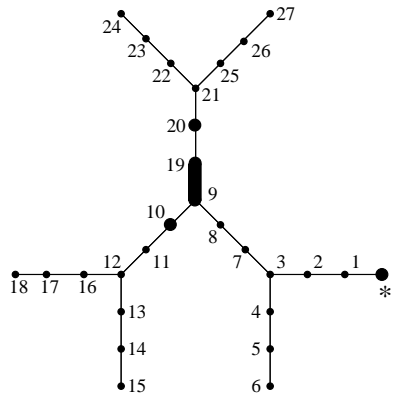}
\caption{An example of a critical $1$-cell\label{fig:critcell}}
\end{center}
\end{figure}

This is the $1$-cell $\{ e_{19}, v_{10}, v_{20}, \ast \}$, by the
numbering of \figref{fig:favtree}.  Notice that, since the order on
vertices depends only on the embedding and the location of $\ast$, it
would be immediately clear that $\tau(e_{19}) < v_{10} < \iota(e_{19})$
even in the absence of an explicit numbering.  It is convenient to
introduce a notation for certain cells (critical cells among them)  which
doesn't refer to a numbering of the vertices.  The notation we introduce
will also be ``independent of subdivision" in a sense which will be
specified.

Fix an alphabet $A, B, C, D, \dots$.  Assign to each essential vertex of
$\Gamma$ a letter from the alphabet.  For this example, we will do this
 so that the smallest essential vertex
is assigned the letter $A$, the next smallest is assigned the letter $B$,
and so on.  Thus, $A = v_{3}$, $B = v_{9}$,
$C= v_{12}$, and $D = v_{21}$.  If $\vec{a} = (a_{1}, \ldots, a_{d(A)-1})$
is a vector with $d(A)-1$ entries, all in $\mathbb{N} \cup \{ 0 \}$, $1
\leq m \leq d(A)-1$, and $a_{m} \neq 0$, we define the notation
$A_{m}[\vec{a}]$ to represent the following subconfiguration.  In the
$m$th direction from $A$, there is a single edge $e$ with $\tau(e) = A$,
and $a_m - 1$ vertices being blocked, stacked up behind $e$.  In any other
direction $i > 0$, $i \neq m$ from $A$, there are $a_i$ blocked vertices
stacked up behind $e$ at $A$. For example, $A_{2}[(1,2)]$ refers to the
collection $\{ v_{4}, e_{7}, v_{8} \}$, and $A_{1}[(3,0)]$ refers to $\{
v_{6}, v_{5}, e_{4} \}$.  We also need a notation for collections of
blocked vertices clustered around $\ast$.  Let $k \ast$ denote a
collection of $k$ vertices blocked at $\ast$.  Thus $1\ast$ refers to $\{
\ast \}$, $2\ast$ refers to $\{ v_{1}, \ast \}$, and so forth.  We combine
these new expressions using additive notation.  Thus, the cell $\{ v_{10},
e_{19}, v_{20}, \ast \}$ would be expressed as $B_{2}[(1,2)] + 1\ast$.  
The critical $2$-cell $\{ v_{13}, e_{16}, v_{10}, e_{19} \}$ would be
expressed as $B_{2}[(1,1)] + C_{2}[(1,1)]$.  If a cell $c$ contains
deleted edges, these can be simply listed, while still using the additive
notation.  For example, one would write $e + k\ast$ to represent a cell
consisting of the deleted edge $e$ together with a collection of $k$
vertices clustered at $\ast$.

We must mention one last convention.  The expression $B_{1}[(4,0)]$ could
refer either to $\{ v_{13}, v_{12}, v_{11}, e_{10} \}$ or $\{ v_{16},
v_{12}, v_{11}, e_{10} \}$.  It is clear that if we subdivide $\Gamma$
further, then this ambiguity disappears.  In arguments using our new
notation, we extend the notion of ``sufficiently subdivided for $n$" in
such a way that all of our new notations are unambiguous for $\Gamma$ when
they refer to collections of $n$ or fewer cells.  It is clearly possible
to do this.

Now our notation is ``subdivision invariant":  if $\Gamma$ is sufficiently
subdivided for $n$, and the expression $A_{m}[\vec{a}] + B_{n}[\vec{b}] +
\ldots + k\ast$ satisfies $a_{1} + \ldots + a_{d(A)-1} + b_{1} + \ldots +
b_{d(B)-1} + \ldots + k \leq n$, then $A_{m}[\vec{a}] + B_{n}[\vec{b}] +
\ldots + k\ast$ specifies a unique cell of $\udng$, and does so no matter
how many times we subdivide.

With these conventions, every $\udng$ has a unique critical $0$-cell,
namely $n\ast$.  The discrete gradient vector field $W$ thus determines a
presentation of $B_{n}\Gamma$ which is unique up to the choices of the
boundary words of critical $2$-cells in $\udng$.

We record for future reference a description of critical cells, in terms
of the new notation:

\begin{proposition} \label{prop:critnotation}

Let $\Gamma$ be a sufficiently subdivided graph, with a maximal subtree
$T$ and basepoint $\ast$. Suppose that $T$ has been embedded in the plane,
so that there is a natural order on the vertices. Assume also that the
endpoint of each deleted edge has degree $1$ in $T$.  With notation as    
above, a cell described by a formal sum
  $$ A_{l_{A}}[\vec{a}] + B_{l_{B}}[\vec{b}] + \ldots + e_{1} + e_{2} +    
  \ldots + k\ast, $$
where each $e_{i}$ is a deleted edge, is critical provided that, for $X=
A,B, \ldots$, some component $x_{j}$ of $\vec{x} = (x_1, x_2, \ldots,
x_{d(X)-1})$ is non-zero, for $j < l_{X}$. Conversely, every critical cell
can be described by such a sum.

If, for some essential vertex $X$, $x_{j}=0$ for all $j < l_{X}$, then  
the above cell is collapsible.\qed

\end{proposition}

We refer to each term of the sum in Proposition \ref{prop:critnotation} 
as a \emph{subconfiguration}.  If the subconfiguratin itself represents a 
critical cell (with fewer strands), we call it \emph{critical}.


\section{Corollaries}


\begin{theorem}\label{thm:numgens}
Let $\Gamma$ be a graph with the maximal tree $T$.  If $T \neq \Gamma$,
then assume that the endpoints of every deleted edge have degree $1$ in
the tree $T$ (and furthermore that $T$ is sufficiently subdivided).  Fix    
the discrete gradient vector field $W$ as in the previous section.  Let    
$D$ be the number of deleted edges.  Then $\mathcal{P}_{W}$ has
  $$D + \sum_{\substack{v \in V(T)\\essential}} \hspace{2mm} \sum_{i =    
  2}^{d(v)-1} \left[\binom{n+d(v)-2}{n-1} - \binom{n+d(v)-i-1}{n-1}\right]$$
generators.
\end{theorem}

\begin{proof}
If $c$ is a critical $1$-cell, then $c$ contains exactly one edge $e$,
which must either be a deleted edge or there is some $v$ in $c$ such that
$e(v) \cap e = \tau(e)$.  In the latter case, $\tau(e)$ is an essential
vertex, and $0 < g(\tau(e), v) < g(\tau(e), \iota(e))$. It follows that $2
\leq g(\tau(e), \iota(e)) \leq d(\tau(e))-1$.

Now let us count the critical $1$-cells $c$.  If the unique edge $e$ in
$c$ is a deleted edge, then $c$ is uniquely determined by $e$, so there
are exactly $D$ critical $1$-cells of this description. If the edge $e$ is
not a deleted edge, then $\tau(e)$ is essential and $2 \leq
g(\tau(e),\iota(e)) \leq d(\tau(e)) - 1$.  Note that since every vertex of
$c$ is necessarily blocked, the critical cell $c$ is determined by the
numbers of vertices of $c$ that are in each of the $d(\tau(e))$ connected
components of $T - e$.  There are $\binom{n+d(\tau(e))-2}{n-1}$ ways to
assign $n-1$ vertices to the $d(\tau(e))$ connected components of $T-e$.     
(This is the number of ways to assign $n-1$ indistinguishable balls to
$d(\tau(e))$ distinguishable boxes.)  Not every such assignment results in
a critical $1$-cell, however. The condition that $0 < g(\tau(e), v) <
g(\tau(e), \iota(e))$, for some $v \in c$,
 won't be satisfied if, for each $v \in c$, either $g(\tau(e), \iota(e)) \leq
g(\tau(e), v) \leq d(\tau(e))-1$ or $g(\tau(e),v) = 0$.  There are
$\binom{d(\tau(e))+ n - g(\tau(e), \iota(e)) - 1}{n-1}$ such ``illegal''
assignments.  Subtracting these from the total, we get \begin{displaymath}
\binom{n + d(\tau(e))-2}{n-1} - \binom{n + d(\tau(e)) -
g(\tau(e),\iota(e)) - 1}{n-1} \end{displaymath} different critical
$1$-cells for a fixed edge $e$.  Letting the edge $e$ of $c$ vary over all
possibilities, we obtain the sum in the statement of the theorem.     
\end{proof}

\begin{corollary}{\rm\cite{G}}\qua
If $\Gamma$ is a \emph{radial} tree -- i.e.\ has exactly one essential
vertex, $v$ -- then $\bng$ is free of rank
  $$\sum_{i=2}^{d(v)-1}
  \left[\binom{n+d(v)-2}{n-1} -   \binom{n+d(v)-i-1}{n-1}\right].$$
\end{corollary}

\begin{proof}
There are no critical cells of dimension greater than $1$ by the
classification of critical cells, since each blocking edge must be at its
own essential vertex.  Thus the presentation $\mathcal{P}_{W}$ has no
relations.  By Theorem \ref{thm:numgens}, the presentation
$\mathcal{P}_{W}$ has the given number of generators.
\end{proof}

\begin{theorem}\label{thm:maxdim}
Let $\Gamma$ be a tree and $c$ a critical cell of $\udng$.  Let
$$k :=  \min \left\{\left\lfloor\frac{n}{2}\right\rfloor, \#\{ v \in \Gamma^{0} 
\mid v \hbox{ is essential}\}\right\}.$$
Then $\dim c \leq k.$ In particular, $\udng$ strong deformation retracts 
on $\left( \udng \right)_{k}'$.
\end{theorem}

\begin{proof}
For every edge $e$ in $c$, there is some vertex $v$ in $c$ such that $e(v)
\cap e = \tau(e)$.  If $e_{1}$ and $e_{2}$ are two edges in $c$, and the
vertices $v_1, v_2 \in c$ satisfy $e(v_i) \cap e_i = \tau(e_i)$, for
$i=1,2$, then certainly $v_1 \neq v_2$, since $e_1 \cap e_2 = \emptyset$.
It follows that there are at least as many vertices as edges in $c$.     
Since the dimension of $c$ is equal to the number of edges in $c$, and the
total number of cells in $c$ is $n$, we have that the dimension of $c$ is
less than or equal to $n/2$.

Since $\tau(e)$ must be an essential vertex of $\Gamma$ for each $e$ in
$c$, and the edges contained in $c$ must be disjoint, we have that the
dimension of $c$ is bounded above by the number of essential vertices of
$\Gamma$.

The final statement now follows from Proposition \ref{prop:corollaries}(2).
\end{proof}

The following result was proven independently for the tree case by Carl 
Mautner, an REU student working under Aaron Abrams \cite{A3}.

\begin{theorem}
Let $\Gamma$ be a sufficiently subdivided graph, and let 
$\chi(\Gamma)$ denote the Euler characteristic of $\Gamma$.  
Then $\udng$ 
strong deformation retracts onto a CW-complex of dimension at most $k$, 
where
  $$k := \min \left\{\left\lfloor\frac{n+1-\chi(\Gamma)}{2}\right\rfloor,
  \#\{v \in \Gamma^{0} \mid v \hbox{ is essential}\}\right\}.$$
\end{theorem}

\begin{proof}
We construct a maximal subtree $T$ of $\Gamma$ whose deleted edges all 
neighbor essential vertices in $\Gamma$.  We note that the existence of a 
\emph{connected} maximal subtree with this property is not clear a priori.

Let $T'$ be any maximal subtree for $\Gamma$.  Let $T$ be a maximal
subtree of $\Gamma$ constructed as follows.  For every deleted edge $e \in
\Gamma - T'$, there are two essential vertices of $\Gamma$ nearest $e$.  
For each deleted edge $e$, let $v_e$ a choice of be such a vertex.  
Construct $T$ by removing from $\Gamma$, for every deleted edge $e$ of
$T'$, the unique edge adjacent to $v_e$ on the simple path from $v_e$ to
$e$ which does not cross any essential vertices.  Then $T$ is connected
and is a maximal tree since $T'$ is.

Any embedding of the tree $T$ induces a discrete gradient vector field $W$
with the property that \emph{every} edge in a critical cell contains an
essential vertex.  Thus, the number of essential vertices bounds the
dimension of the cells of $\udng$ that are critical with respect to $W$.

Let $D$ be the number of deleted edges with respect to $T$.  Thus $D = 1 -
\chi(\Gamma)$.  Since a critical subconfiguration involves either one
strand on a deleted edge or at leasttwo strands about an essential vertex,
the dimension of any critical cell is also bounded by $D +
\lfloor\frac{n-D}{2}\rfloor = \lfloor\frac{n+D}{2}\rfloor$.

The theorem now follows from Proposition \ref{prop:corollaries}(2).
\end{proof}

\begin{figure}[ht!]\anchor{fig:k_4}
\begin{center}
\includegraphics{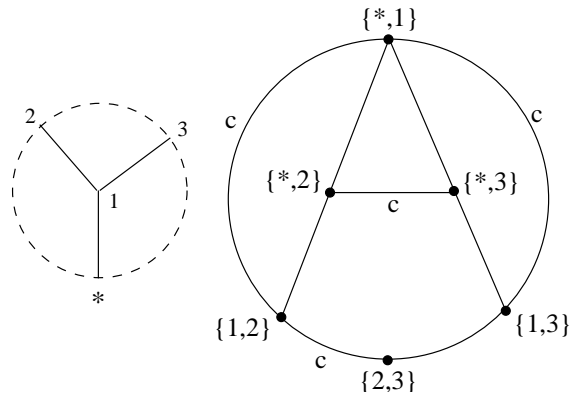}
\caption{On the left, we have $K_4$ with the choice of maximal tree
depicted.  On the right, we have the complex $X_1'$ onto which
$\mathcal{UD}^2K_4$ deformation retracts.  The critical $1$-cells are
indicated with a lower case ``c"; the critical $0$-cell 
is $\{*,1\}$.\label{fig:k_4}}
\end{center}
\end{figure}

Note that this bound is not sharp, as two deleted edges for $T$ adjacent    
to the same essential vertex may not both be involved in a critical    
configuration.  The real thing to count (which varies by choice of $T$) is    
the number of essential vertices in $\Gamma$ which touch deleted    
edges, not deleted edges themselves.  

\begin{example}
Consider $\mathcal{UD}^{2}K_{4}$, where $K_{4}$ is the complete graph on
four vertices.  We choose a radial tree as our maximal tree.  It is not
difficult to check that there are no critical $2$-cells in
$\mathcal{UD}^{2}K_{4}$ with respect to this choice of maximal tree, so
the subcomplex $X_{1}' \subseteq \mathcal{UD}^{2}K_{4}$ is a strong
deformation retract of $\mathcal{UD}^{2}K_{4}$ (see \figref{fig:k_4}).
\end{example}


\section{Presentations of tree braid groups}


Let $v$ and $v'$ be vertices of the graph $\Gamma$.  Let $(v, v') = \{ v''
\in V(\Gamma) \mid v < v'' < v' \}$; in other words, $(v,v')$ denotes the
set of vertices ``between'' $v$ and $v'$ in the ordering on vertices.     
The following lemma is extremely useful in computing presentations of
graph braid groups.

\begin{lemma}[Redundant 1-cells lemma]
Let $c$ be a redundant $1$-cell in\break $\udng$, where $\Gamma$ is an arbitrary
finite graph. Suppose $c = \{c_1, \dots, c_{n-2}, v,e\}$, where $v$
is the smallest unblocked vertex of $c$ and $e$ is the unique edge of $c$.
Let $c' := \{c_1, \dots, c_{n-2},\tau(e(v)), e\}$.  If $(\tau(e(v)),v)
\cap \{ c_1, c_2, \ldots, c_{n-2}, \iota(e), \tau(e) \} = \emptyset$, then
$c \dot{\rightarrow} c'$.  

Here $\dot{\rightarrow}$ refers to a sequence of moves over the complete
rewrite system $\mathcal{MP}_{W,T}$ (see Proposition \ref{cmplt} and the 
discussion preceding).  In particular, $c$ and $c'$ are words (of length $1$)
in the alphabet of oriented $1$-cells.   
\end{lemma}

\begin{proof}
Let $c_{\iota} := \{ c_1, \ldots, c_{n-2}, e(v), \iota(e) \}$ and
$c_{\tau} := \{ c_1, \ldots, c_{n-2}, e(v), \tau(e) \}$.  If we apply $W$
to $c$, then we get the following $2$-cell, where the arrows indicate
orientation:
\[
\xymatrix@R=20pt@C=30pt{
  *=0{} \ar@{-}|{\object@{>}}[rr]^{c}
        \ar@{-}|{\object@{>}}[dd]_{c_{\iota}} & &
  *=0{} \ar@{-}|{\object@{>}}[dd]^{c_{\tau}} \\ \\
  *=0{} \ar@{-}|{\object@{>}}[rr]_{c'} & &
  *=0{} } \]
It follows that $c \rightarrow c_{\iota}c'\overline{c_{\tau}}$; we need to
show that $c_{\iota} \dot{\rightarrow} 1$ and $c_{\tau} \dot{\rightarrow}
1$.

Consider the following condition on $1$-cells $\hat{c}$ of $\udng$:

\begin{enumerate}
\item[(1)]  If $v$ is a vertex and $v \in \hat{c}$, then      
$v \not \in (\tau(\hat{e}), \iota(\hat{e}))$, where $\hat{e}      
\subseteq T$ is the unique edge of $\hat{c}$.
\end{enumerate}

We claim: first, that if a 1-cell $\hat{c}$ satisfies $(1)$, then so does
any 1-cell in a $W$-path starting with $\hat{c}$; second, that a 1-cell
$\hat{c}$ satisfying $(1)$ cannot be critical.  Since $c_{\iota}$ and
$c_{\tau}$ both satisfy $(1)$ by hypothesis, this will prove the lemma by
induction on the rank of $c_{\iota}$ and $c_{\tau}$.

We begin with the latter claim.  Suppose $\hat{c}$ satisfies $(1)$. If
$\hat{c}$ is critical, then there is $v \in \hat{c}$ such that $e(v) \cap
\hat{e} = \tau(\hat{e})$, and $0 < g(\tau(\hat{e}),v) < g(\tau(\hat{e}),
\iota(\hat{e}))$.  But then $\tau(\hat{e}) < v < \iota(\hat{e})$, a
contradiction.

Now for the first claim.  If $\hat{c}$ is collapsible, then any $W$-path
starting with $\hat{c}$ consists of $\hat{c}$ alone.  We may therefore
assume that $\hat{c}$ is redundant.  Since $\hat{c}$ is redundant, there
is $v \in \hat{c}$ such that $v < \iota(\hat{e})$ and $v$ is unblocked.     
We may choose $v$ to be the minimal vertex with this property.  By
assumption, $v < \tau(\hat{e})$ as well.  The two edges in $W(\hat{c})$,
namely $e(v)$ and $\hat{e}$, therefore have the property that $\tau(e(v))
< \iota(e(v)) < \tau(\hat{e}) < \iota(\hat{e})$. Moreover, there is no
vertex in $W(\hat{c}) \cap (\tau(e(v)),\iota(e(v)))$, for a minimal such
vertex $v'$ would be an unblocked vertex in $\hat{c}$ satisfying $v' < v$,
which is impossible.  The immediate successor $\hat{c}_{1}$ of $\hat{c}$
in a $W$-path is obtained by replacing an edge of $W(\hat{c})$ with either
its initial or terminal vertex.  It thus follows that $\hat{c}_{1}$
satisfies $(1)$.
\end{proof}

In order to prove our theorem on the presentations of tree braid groups,
we need to introduce new notation.  Let $\dot{A}[\vec{a}]$ denote the
collection of vertices consisting of $A$ itself together with $a_{i}$
vertices arranged consecutively in the $i$th direction from $A$, so that
every vertex in the collection is blocked except for $A$.  Thus, to use
the tree from \figref{fig:favtree}, the expression $\dot{B}[(2,1)]$
would refer to the collection $\{ v_{9}, v_{10}, v_{11}, v_{19} \}$.  Let
$A[\vec{a}]$ denote the same collection as does $\dot{A}[\vec{a}]$, but
excluding the vertex $A$ itself.  Thus, $B[(2,1)]$ refers to $\{ v_{10},
v_{11}, v_{19} \}$ in our favorite tree. We combine these new notations
with the old ones (namely, $A_{m}[\vec{a}]$ and $k\ast$) additively, as
before.

A word about notation:  we will need to describe the boundaries of
critical $2$-cells in $\udng$, and these boundaries sometimes consist of
certain $1$-cells which cannot be expressed in terms of our notation.     
For instance, consider $B_{2}[(1,2)] + D_{2}[(1,1)]$ in relation to 
\figref{fig:favtree}. This is the $2$-cell $c=\{ v_{10}, e_{19}, v_{20},
v_{22}, e_{25} \}$.  The $1$-cells forming the boundary are obtained by
replacing the edges $e_{i} \in c$, for $i = 19$ or $25$, with
$\iota(e_{i})$ or $\tau(e_{i})$.  Thus one of the $1$-cells on the
boundary of $c$ is $\{ v_{10}, v_{9}, v_{20}, v_{22}, e_{25} \}$.  Note
that this configuration isn't covered by our notation, since the
inessential vertex $v_{20}$ is unblocked.

Rather than introduce more notation to deal with these extra
configurations, we introduce the idea of a \emph{slide}.  A $1$-cell $c'$
is obtained from $c$ by a slide (or $c$ \emph{slides to} $c'$) if there is
some unblocked vertex $v \in c$ such that the endpoints, say $v$ and $v'$,
of $e(v)$ are consecutive in the order on vertices, and $c'$ is obtained
from $c$ by replacing $v$ with $v'$.  For example, $\{ v_{10}, v_{9},
v_{19}, v_{22}, e_{25} \}$ (which is $\dot{B}[(1,1)] + D_{2}[(1,1)]$) is
obtained from $\{ v_{10}, v_{9}, v_{20}, v_{22}, e_{25} \}$ by a slide.

If $c$ slides to $c'$, then we may use them interchangeably in our
calculations.  For if $c = \{ c_1, \ldots, c_{n-2}, v, e \}$ and $c' = \{
c_1, \ldots, c_{n-2}, v', e \}$, where $v' = \tau(e(v))$, then $c$ and
$c'$ are parallel sides of the square $\{ c_{1}, \ldots, c_{n-2}, e(v), e
\}$, and the other sides are $c_{\iota} = \{ c_{1}, \ldots, c_{n-2}, e(v),
\iota(e) \}$ and $c_{\tau} = \{ c_{1}, \ldots, c_{n-2},$
$ e(v), \tau(e) \}$.
It is clear that $c_{\iota}$ and $c_{\tau}$ both satisfy the condition
$(1)$ from the proof of the redundant 1-cells lemma, so $c_{\iota}
\dot{\rightarrow} 1$ and $c_{\tau} \dot{\rightarrow} 1$.  It then follows
that the oriented 1-cells $c$ and $c'$ both represent the same element in
the usual edge-path presentation of $\pi_{1}(\udng)$.  In fact, more is
true: if we add all relations corresponding to slides $(c,c')$ to the
monoid presentation $\mathcal{MP}_{W,T}$, the associated string rewriting
system is still complete, and has the same reduced objects.  We leave the
verification of this fact as an exercise.  (Note that the square $\{
c_{1}, \ldots, c_{n-2}, e(v), e \}$ may be redundant, so the proof is not
entirely trivial.)  In our calculations, we will use slides without
further notice.

Let $\delta_{m}$ denote a vector such that the $m$th component is $1$ and
every other component is $0$.  The length of $\delta_{m}$ will be clear
from the context.  If $\vec{v}$ is a vector having entries in the set of
non-negative integers, let $\vec{v} - 1$ be the vector obtained from
$\vec{v}$ by subtracting $1$ from the first non-zero entry.  This last
notation must be used carefully to avoid ambiguity -- note for instance
that $\delta_{1} + (\delta_{2} - 1) \neq (\delta_{1} + \delta_{2}) - 1$.     
If $\vec{v}$ is any vector, let $|\vec{v}|$ denote the sum of the entries
of $\vec{v}$.

\begin{lemma}
Let $A$ and $B$ be essential vertices of $T$, a maximal tree in $\Gamma$.
\begin{enumerate}
\item Suppose that $A \wedge B = C$ where $C$ is an essential vertex    
distinct from both $A$ and $B$. Let $g(C,A)=i$ and $g(C,B)=j$, where $i <    
j$.  Then

\begin{enumerate}[\rm(a)]
\item $\left(A[\vec{a}] + B_{l}[\vec{b}] + p\ast\right) \dot{\rightarrow}    
\left(B_{l}[\vec{b}] + (p +  |\vec{a}|)\ast\right)$,

\item $\left(\dot{A}[\vec{a}] + B_{l}[\vec{b}] + p\ast\right)    
\dot{\rightarrow} \left(B_{l}[\vec{b}] + (p + 1 + |\vec{a}|)\ast\right)$,

\item $\left(A_{k}[\vec{a}] + B[\vec{b}] + p\ast\right) \dot{\rightarrow}
w_{1}\left(A_{k}[\vec{a}] + (p+  |\vec{b}|)\ast\right)w_{1}^{-1}$, and

\item $\left(A_{k}[\vec{a}] + \dot{B}[\vec{b}] + p\ast\right)    
\dot{\rightarrow} w_{2}\left(A_{k}[\vec{a}] +    
(1+p+|\vec{b}|)\ast\right)w_{2}^{-1}$,
\end{enumerate}
where    
  $$w_{1} = \prod_{\alpha = 0}^{|\vec{b}|-1}     
  \left(C_{j}[|\vec{a}|\delta_{i}  + (|\vec{b}| - \alpha)\delta_{j}] + (p    
  + \alpha)\ast\right),$$
and $w_{2}$ is the same as $w_{1}$, but with $|\vec{b}|+1$ in place of    
$|\vec{b}|$.

\item Suppose that $A \wedge B = A$ and $g(A,B) = i$.  Then
\begin{enumerate}[\rm(a)]
\item $\left(A_{k}[\vec{a}] + B[\vec{b}] + p\ast\right) \dot{\rightarrow}    
\left(A_{k}[\vec{a} + |\vec{b}|\delta_{i}] + p\ast\right)$,

\item $\left(A_{k}[\vec{a}] + \dot{B}[\vec{b}] + p\ast\right)
\dot{\rightarrow} \left(A_{k}[\vec{a} + (1 + |\vec{b}|)\delta_{i}] +    
p\ast\right)$,

\item $\left(A[\vec{a}] + B_{l}[\vec{b}] + p\ast\right) \dot{\rightarrow}      
w_{3}\left(B_{l}[\vec{b}] + (p+|\vec{a}|)\ast\right)w_{3}^{-1}$, and

\item $\left(\dot{A}[\vec{a}] + B_{l}[\vec{b}] + p\ast\right)    
\dot{\rightarrow} \left(A[\vec{a}] + B_{l}[\vec{b}] + (p+1)\ast\right)$,
\end{enumerate}
where
  $$w_{3} = \prod_{\alpha=0}^{|\vec{a}|-1}    
  \left(A_{\beta}[|\vec{b}|\delta_{i} + (\vec{a}- \alpha)] + (p+    
  \alpha)\ast\right),$$
and $\beta$ is the smallest coordinate of $\vec{a} - \alpha$ that is    
non-zero.  Here a factor in $w_3$ is considered trivial if $\beta \leq i$.
\end{enumerate}
\end{lemma}

\begin{proof}$\phantom{99}$

(1a)\qua Under the given assumptions, the smallest vertex of the    
subconfiguration $A[\vec{a}]$ may be moved until it is blocked at $\ast$, 
by repeated applications of the redundant 1-cells lemma.  That is,
  $$ \left( A[\vec{a}] + B_{l}[\vec{b}] + p\ast \right) \dot{\rightarrow} 
\left( A[\vec{a}-1] +    
  B_{l}[\vec{b}] + (p+1)\ast \right). $$
After repeated applications of the above identity, we eventually arrive    
at the statement (a).

(1b)\qua This is similar to (a).

(1c)\qua We begin by applying the redundant 1-cells lemma to the
smallest vertex of the subconfiguration $B[\vec{b}]$.  This smallest 
vertex can be moved freely, until it occupies the place adjacent to $C$, 
and lying in the $j$th direction from $C$.  That is,
  $$ \left( A_{k}[\vec{a}] + B[\vec{b}] + p\ast \right) \dot{\rightarrow}
  \left( A_{k}[\vec{a}]+B[\vec{b}-1]+C[\delta_{j}]+p\ast \right).$$
At this point, the redundant 1-cells lemma no longer applies, since all of    
the vertices in the subconfiguration $A_{k}[\vec{a}]$ lie between $C$ and 
the vertex (say $v$) which is adjacent to $C$ and which lies in the $j$th    
direction.  We must therefore appeal to the definition of $W$:         
\[
\xymatrix@R=35pt@C=40pt{
  *=0{} \ar@{-}|{\object@{>}}[rr]^{A_{k}[\vec{a}] + B[\vec{b}-1] + C[\delta_{j}] + p\ast}
        \ar@{-}|{\object@{>}}[dd]_{A[\vec{a}] + B[\vec{b}-1] + C_{j}[\delta_{j}] + p\ast} & &
  *=0{} \ar@{-}|{\object@{>}}[dd]^{\dot{A}[\vec{a}-\delta_{k}] + B[\vec{b}-1] + C_{j}[\delta_{j}] + p\ast} \\ \\
  *=0{} \ar@{-}|{\object@{>}}[rr]_{A_{k}[\vec{a}] + B[\vec{b}-1] + \dot{C} + p\ast} & &
  *=0{}
}
\]
The $2$-cell depicted above is $A_{k}[\vec{a}] + B[\vec{b}-1] +
C_{j}[\delta_{j}] + p\ast$, the image under $W$ of the $1$-cell
$A_{k}[\vec{a}] + B[\vec{b}-1] + C[\delta_{j}] + p\ast$ (located at the
top).  (Note:  the label of the ``source" vertex on the upper left can be
computed by replacing each of the edges in the $2$-cell $A_{k}[\vec{a}] +
B[\vec{b}-1] + C_{j}[\delta_{j}] + p\ast$ with its initial vertex.  The
``sink" vertex on the bottom right is obtained from the same $2$-cell by
replacing each edge with its terminal vertex.  The other two vertices are
determined by replacing one of the edges with its initial vertex, and the
other with its terminal vertex.  Furthermore, if we specify the identity 
of any one of the 1-cells on the boundary of the 2-cell, then the labels 
of the other three 1-cells are uniquely determined.)

Now consider the 1-cell $A[\vec{a}]+B[\vec{b}-1]+C_{j}[\delta_{j}]+p\ast$.  
The redundant 1-cells lemma applies to the vertices in the
subconfiguration $A[\vec{a}]$.  These may move until they are blocked at
$C$.  Since there are $|\vec{a}|$ vertices in $A[\vec{a}]$, and each lies
in the direction $i$ from $C$, we have
  $$ \left( A[\vec{a}]+B[\vec{b}-1]+C_{j}[\delta_{j}]+p\ast \right) 
\dot{\rightarrow}    
  \left( B[\vec{b}-1] + C_{j}[\delta_{j} + |\vec{a}|\delta_{i}] + p\ast \right).$$
If we let the vertices in the subconfiguration $B[\vec{b}-1]$ move, then 
by the redundant 1-cells lemma, these vertices will flow until they are 
blocked at $C$.  We get
  $$ \left( B[\vec{b}-1]+C_{j}[\delta_{j} + |\vec{a}|\delta_{i}] + p\ast \right)    
  \dot{\rightarrow} \left( C_{j}[|\vec{b}|\delta_{j} + |\vec{a}|\delta_{i}] +    
  p\ast \right).$$
This last 1-cell is critical.  The 1-cell on the right side of the square
pictured above flows to the same 1-cell, by similar reasoning.

Finally, the 1-cell $A_{k}[\vec{a}] + B[\vec{b}-1]+ \dot{C} + p\ast$ flows    
to $A_{k}[\vec{a}] + B[\vec{b}-1] + (p+1)\ast$ by the redundant 1-cells    
lemma.  It follows that
  $$ \left( A_{k}[\vec{a}] + B[\vec{b}] + p\ast \right) \dot{\rightarrow}
   w \left( A_{k}[\vec{a}]+B[\vec{b}-1]+(p+1)\ast \right) w^{-1},$$
where $w = (C_{j}[|\vec{a}|\delta_{i} + |\vec{b}|\delta_{j}] + p\ast)$.
Part (c) now follows by repeatedly applying the above identity.

(1d)\qua This is similar to (c).

(2a)\qua
In this case, the vertices in the subconfiguration $B[\vec{b}]$, beginning
with the smallest, can flow by the redundant 1-cells lemma until they are
blocked at the essential vertex $A$.  Since the vertices in $B[\vec{b}]$
lie in the direction $i$ from $A$, when these vertices are blocked, the
resulting configuration is $A_{k}[\vec{a}+|\vec{b}|\delta_{i}] + p\ast$.     
This proves (a).
     
(2b)\qua This is similar to (a).
     
(2c)\qua We begin by applying $W$ to the configuration $A[\vec{a}] +
B_{l}[\vec{b}] + p\ast$.  The result is $A_{\beta}[\vec{a}] +
B_{l}[\vec{b}] + p\ast$, where $\beta$ is the smallest subscript for which
$a_{\beta}$ is non-zero (here $\vec{a} = (a_1, a_2, \ldots, a_{d(A)-1})$):
\[
\xymatrix@R=35pt@C=40pt{
  *=0{} \ar@{-}|{\object@{>}}[rr]^{A[\vec{a}] + B_{l}[\vec{b}] + p\ast}
        \ar@{-}|{\object@{>}}[dd]_{A_{\beta}[\vec{a}] + B[\vec{b}] + p\ast} & &
  *=0{} \ar@{-}|{\object@{>}}[dd]^{A_{\beta}[\vec{a}] + \dot{B}[\vec{b}-\delta_{l}] + p\ast} \\ \\
  *=0{} \ar@{-}|{\object@{>}}[rr]_{\dot{A}[\vec{a}-1] + B_{l}[\vec{b}] + p\ast} & &
  *=0{}
}
\]
Now it is either the case that both of the vertical 1-cells are 
collapsible (if $\beta \leq i$), or these vertical 1-cells both flow to     
  $$ \left( A_{\beta}[\vec{a} + |\vec{b}|\delta_{i}] + p\ast \right) . $$
The bottom 1-cell flows to     
  $$ \left( A[\vec{a}-1]+ B_{l}[\vec{b}] + (p+1)\ast \right) .$$
Thus, $\left( A[\vec{a}] + B_{l}[\vec{b}] + p\ast \right) $ flows to     
  $$ w\left( A[\vec{a}-1]+B_{l}[\vec{b}]+(p+1)\ast \right) w^{-1},$$
where $w = \left( A_{\beta}[\vec{a} + |\vec{b}|\delta_{i}] + p\ast \right)$ 
if $\beta >    
i$, and $w=1$ if $\beta \leq i$.  The statement of 2(c) now follows by    
repeated application of the above identity.

(2d)\qua This is straightforward.
\end{proof}

\begin{theorem} \label{thm:treepres}
Let $\Gamma$ be a sufficiently subdivided tree with a chosen basepoint
$\ast$.  Suppose that an embedding of $\Gamma$ in the plane is given, so
that there is an induced order on the vertices.  Then the braid group 
$\bng$ is generated by the collection of critical $1$-cells, and the set 
of relations consists of the reduced forms of the boundary words $w(c)$, 
where $c$ is any critical $2$-cell.  For $c = \left(A_{k}[\vec{a}] + 
B_{l}[\vec{b}] + p\ast\right)$ a critical $2$-cell, the reduced form of 
the boundary word is as follows:

\begin{enumerate}
\item If $A \wedge B = C$ with $C \neq A,B$, $g(C,A) = i$, $g(C,B) = j$, 
and $i<j$, then a reduced form of the boundary word for $c$ is
  $$ \left[\left(B_{l}[\vec{b}] + (p+|\vec{a}|)\ast\right),    
  w_{1}\left(A_{k}[\vec{a}] + (p+|\vec{b}|)\ast\right)w_{1}^{-1} \right],$$
where $w_{1}$ is as in the previous lemma.
\item If $A \wedge B = A$, and $g(A,B) = i$, then a reduced form of the    
boundary word for $c$ is
  $$ \left[ w_{3}^{-1}\left(A_{k}[\vec{a} + |\vec{b}|\delta_{i}] +    
  p\ast\right)w_{3}', \left(B_{l}[\vec{b}] + (p+    
  |\vec{a}|)\ast\right)^{-1} \right],$$
where $w_{3}$ is as in the previous lemma, and $w_{3}'$ is the same as    
$w_{3}$, but with $\vec{a} - \delta_{k}$ in place of $\vec{a}$ and $p+1$    
in place of $p$.
\end{enumerate}
\end{theorem}

\begin{proof}
It follows from Theorem \ref{thm:Morsepresentation} that $\bng$ is
generated by the critical $1$-cells, and that the relations are the
reduced forms of the boundary words of the critical $2$-cells.

Consider the critical $2$-cell $\left( A_{k}[\vec{a}] + B_{l}[\vec{b}] +    
p\ast \right)$:
\[
\xymatrix@R=35pt@C=40pt{
  *=0{} \ar@{-}|{\object@{>}}[rr]^{A[\vec{a}] + B_{l}[\vec{b}] + p\ast}
        \ar@{-}|{\object@{>}}[dd]_{A_{k}[\vec{a}] + B[\vec{b}] + p\ast} & &
  *=0{} \ar@{-}|{\object@{>}}[dd]^{A_{k}[\vec{a}] + \dot{B}[\vec{b}-\delta_{l}] + p\ast} \\ \\
  *=0{} \ar@{-}|{\object@{>}}[rr]_{\dot{A}[\vec{a}-\delta_{k}] + B_{l}[\vec{b}] + p\ast} & &
  *=0{}
}
\]
the rest of the statement of the theorem follows by applying the previous    
lemma to each side.
\end{proof}

\begin{example}\label{exam:favtree4}
Consider the example of $B_{4}\Gamma$, where $\Gamma$ is the 
tree in \figref{fig:favtree}.  This tree is sufficiently subdivided as 
written.  By Theorem \ref{thm:classification} the critical $1$-cells are:
  $$\left( X_{2}[\vec{v}] + (4 - |\vec{v}|)\ast \right),$$
where $X$ is one of the essential vertices $A$, $B$, $C$, or $D$, and    
$\vec{v}$ is $(1,1)$, $(1,2)$, $(1,3)$, $(2,1)$, $(2,2)$, or $(3,1)$.     
Thus, there are a total of $24$ critical $1$-cells. The critical $2$-cells    
all have the form
  $$\left( X_{2}[(1,1)] + Y_{2}[(1,1)] \right),$$
where $X$ and $Y$ are distinct essential vertices chosen from the set 
$\{A, B, C, D\}$.  Thus there are a total of $6$ critical $2$-cells.

We now describe the resulting relations.  We begin with the critical
$2$-cell $\left(A_{2}[(1,1)] + B_{2}[(1,1)] \right)$. The boundary word,
by (2) of the previous theorem with $i=k=l=2$ and $p=0$, is
  $$ \left[ w_{3}^{-1} \left(A_{2}[(1,3)] \right) w_{3}', \left(
  B_{2}[(1,1)] + 2\ast \right) \right]. $$
The word $w_{3}$ in the present case is
  $$ \prod_{\alpha = 0}^{1} \left(A_{\beta}[(0,2) + ((1,1) - \alpha)] +    
  \alpha \ast \right).$$
Computing, and using the fact that $\beta$ is the first non-zero entry of    
$(1,1) - \alpha$, we get
  $$w_{3} = \left( A_{1}[(1,3)] \right) \left( A_{2}[(0,3)] + 1\ast    
  \right)= 1,$$
since both of the cells in the above expression are collapsible. The 
word $w_{3}'$ in the present case is     
  $$ \prod_{\alpha = 0}^{0} \left( A_{\beta}[(0,2) + \left( (1,0) - \alpha    
  \right)] + (\alpha + 1) \ast \right) = \left( A_{1}[(1,2)] + 1\ast    
  \right) = 1,$$
since the given cell is again collapsible.  Summarizing, we get
  $$ \left[ \left( X_{2}[(1,3)] \right), \left( Y_{2}[(1,1)] + 2\ast    
  \right) \right],$$
where $(X,Y) = (A,B)$.
A completely analogous computation shows that a boundary word for     
$\left( X_{2}[(1,1)] + Y_{2}[(1,1)] \right)$ is given by the same    
expression,
where $(X,Y) \in \{ (A,B), (A,C), (A,D), (B,D) \}$.

Now consider the critical $2$-cell $\left(B_{2}[(1,1)] + C_{2}[(1,1)]
\right)$.  Part (2) of the previous theorem applies again, with $B$
playing the role of $A$ in the theorem, $C$ the role of $B$, $i=1$,
$k=l=2$, and $p=0$:
  $$ \left[ w_{3}^{-1}\left( B_{2}[(3,1)] \right) w_{3}', \left(    
  C_{2}[(1,1)] + 2\ast \right) \right].$$
In the present case, we get the following expression for $w_{3}$:
\begin{eqnarray*}
w_{3} & = & \prod_{\alpha=0}^{1} \left( B_{\beta}[(2,0) + ((1,1) - \alpha)] + \alpha \ast \right) \\
      & = & \left(B_{1}[(3,1)]\right) \left(B_{2}[(2,1)] + 1\ast \right) \\
      & = & \left(B_{2}[(2,1)] + 1\ast \right).
\end{eqnarray*}
We get the following expression for $w_{3}'$:
\begin{eqnarray*}     
w_{3}' & = & \prod_{\alpha=0}^{0} \left( B_{\beta}[(2,0) + ((1,0) - \alpha)] + (\alpha + 1) \ast \right)\\
       & = & \left( B_{1}[(3,0)] + 1\ast \right) \\
       & = & 1.
\end{eqnarray*}
We arrive at the following boundary word for $\left( B_{2}[(1,1)] +    
C_{2}[(1,1)] \right)$:
  $$ \left[ \left( B_{2}[(2,1)] + 1\ast \right)^{-1}\left( B_{2}[(3,1)]    
  \right), \left( C_{2}[(1,1)] + 2\ast \right) \right].$$     
Finally, we consider the critical $2$-cell $\left( C_{2}[(1,1)] +    
D_{2}[(1,1)] \right)$.  By part (1) of the previous theorem, with $C$    
playing the role of $A$, $D$ playing the role of $B$, and $B$ playing the    
role of $C$, $i=1$, $j=k=l=2$, and $p=0$, a boundary word has the form
  $$ \left[ \left( D_{2}[(1,1)] + 2\ast \right), w_{1}\left( C_{2}[(1,1)]    
  + 2\ast \right)w_{1}^{-1} \right],$$
with
\begin{eqnarray*}
w_{1} & = & \prod_{\alpha=0}^{1} \left( B_{2}[(2,2-\alpha)] + \alpha \ast \right) \\
& = & \left( B_{2}[(2,2)] \right) \left( B_{2}[(2,1)] + 1\ast \right).
\end{eqnarray*}
We summarize these calculations in a figure, which also illustrates a    
useful chalkboard  notation for critical $1$-cells (see 
\figref{fig:relators}).    

\begin{figure}[ht!]\anchor{fig:relators}
\begin{center}
\psfrag{,}{\large,}
\psfraga <-2pt, 0pt> {A}{\footnotesize$A$}
\psfraga <-2pt, 0pt> {B}{\footnotesize$B$}
\psfraga <-2pt, 0pt> {C}{\footnotesize$C$}
\psfraga <-2pt, 0pt> {D}{\footnotesize$D$}
\includegraphics{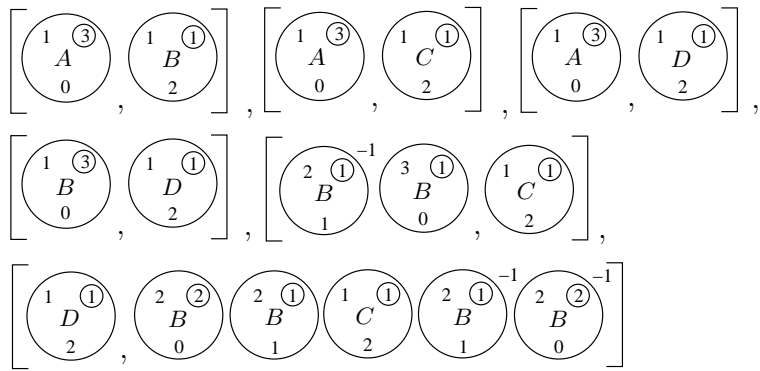}
\caption{The relations of $B_{4}\Gamma$ in shorthand 
notation\label{fig:relators}}
\end{center}
\end{figure}

For example, the circled letter $C$ with a $1$, a circled $1$, and a $2$
(reading in the clockwise direction from the upper left) represents the
element $C_{2}[(1,1)] + 2\ast$: the first entry in the upper left refers
to the first entry of the vector $(1,1)$, the circled entry in the upper
right refers to the second entry of $(1,1)$ and indicates the location of
the edge, and the $2$ on the bottom refers to the number of vertices
clustered near $\ast$.
\end{example}

A group $G$ is said to be a \emph{right-angled Artin group} provided there
is a presentation $\mathcal{P} = \langle \Sigma \mid \mathcal{R} \rangle$
such that every relation in $\mathcal{R}$ has the form $[a,b]$, where $a,
b \in \Sigma$.  It is common to describe a right-angled Artin group
presentation by a graph $\Gamma_{\mathcal{P}}$, where the vertices of
$\Gamma_{\mathcal{P}}$ are in one-to-one correspondence with elements of
$\Sigma$, and there is an edge connecting two vertices $a,b \in \Sigma$ if
and only if $[a,b] \in \mathcal{R}$.

A finitely presented group $G$ is \emph{coherent} if every finitely
generated subgroup of $G$ is also finitely presentable.

\begin{example}\label{exam:Htree} Consider the case in which $\Gamma$ is a
tree homeomorphic to the capital letter ``H".  Let the basepoint $\ast$ be
the vertex on the bottom left, let $A$ be the essential vertex on the
left, and let $B$ be the essential vertex on the right.

\begin{figure}[ht!]\anchor{fig:rapres}
\begin{center}
\psfraga <-1.5pt,0pt> {A}{\tiny$A$}
\psfraga <-1.5pt,0pt> {B}{\tiny$B$}
\includegraphics[width=.9\hsize]{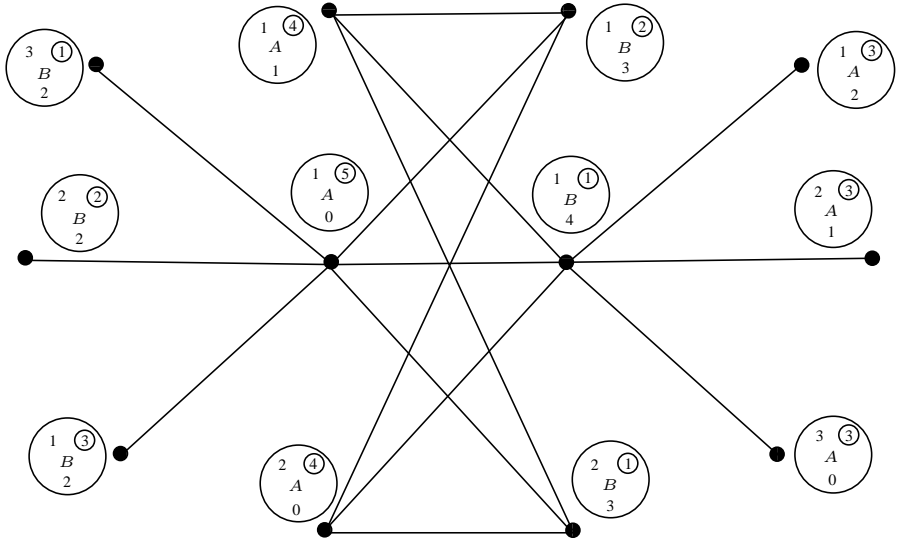}   
\caption{Part of the presentation for $B_{6}H$ \label{fig:rapres}}
\end{center}
\end{figure}

We claim that $\bng$ is a right-angled Artin group, for any $n$.  (See
also Connolly and Doig \cite{CD}, who prove, more generally, that the
braid group of any linear tree is a right-angled Artin group.  A tree $T$
is \emph{linear} if every essential vertex of $T$ lies along a single
embedded arc.)  No matter what $n$ is, part (2) of Theorem
\ref{thm:treepres} applies, and we must compute the values of the words
$w_{3}$ and $w_{3}'$.  In fact, we prove that both of these words will
necessarily be trivial.  Consider $$ w_{3} =
\prod_{\alpha=0}^{|\vec{a}|-1} \left(A_{\beta}\left[ |\vec{b}|\delta_{i} +
\left( \vec{a} - \alpha \right) \right] + \left( p + \alpha \right)\ast
\right).$$ Here $i = g(A,B) = 2$ and $\beta$ is the subscript of the
smallest non-zero entry of $\left( \vec{a} - \alpha \right)$.  The crucial
observation is that $\beta$ is (in the present case) also the smallest
non-zero entry of $|\vec{b}|\delta_{i} + \left( \vec{a} - \alpha \right)$,
so, by Proposition \ref{prop:critnotation}, no term in the above product
is a critical 1-cell (in fact, it is easy to see that each is
collapsible).  The triviality of $w_3$ follows.  Essentially the same
reasoning holds for $w_{3}'$.

If $n=6$, a routine but lengthy calculation shows that
$B_{n}\Gamma$ is the free product of a free group of rank $18$
with the right-angled Artin group in \figref{fig:rapres}.

The copy of $K_{3,3}$ in this graph represents a subgroup isomorphic to
 $F_{3} \times F_{3}$, the direct product
of the free group of rank $3$ with itself.  Since $F_{3} \times F_{3}$ is not coherent (see, for instance,
\cite{BB}) and  
has an unsolvable generalized word
problem \cite{unsol}, there are tree braid groups which are not coherent and have an unsolvable generalized
word problem.
     
\end{example}   

\def\polhk#1{\setbox0=\hbox{#1}{\ooalign{\hidewidth
  \lower1.5ex\hbox{`}\hidewidth\crcr\unhbox0}}}

\Addresses\recd

\end{document}

%% file: agtout.tex
\def\ifplaintex{\expandafter\ifx\csname documentclass\endcsname\relax}

\def\gtp{{\mathsurround=0pt\it $\cal G\mskip-2mu$eometry \&\ 
$\cal T\!\!$opology $\cal P\!$ublications}}  

\def\recd{{\small Received:\qua\receiveddate\ifx\reviseddate\relax
\else\qquad Revised:\qua\reviseddate\fi\par}} 

\def\lognumber#1{\def\thelognumber{#1}}
\def\volumenumber#1{\def\thevolumenumber{#1}}
\def\volumeyear#1{\def\thevolumeyear{#1}}
\def\papernumber#1{\def\thepapernumber{#1}}
\def\pagenumbers#1#2{\def\startpage{#1}\def\finishpage{#2}}
\def\published#1{\def\publishdate{#1}}

\def\received#1{\def\receiveddate{#1}}

\def\accepted#1{\def\accepteddate{#1}}

\def\asciiemail#1{\def\theasciiemail{#1}}
\def\asciiurl#1{\def\theasciiurl{#1}}

\long\def\asciiabstract#1{\long\def\theasciiabstract{#1}}

\let\\\par\let\thelognumber\relax\let\thevolumenumber\relax
\let\thepapernumber\relax\let\thevolumeyear\relax\let\startpage\relax
\let\finishpage\relax\let\publishdate\relax\let\receiveddate\relax
\let\reviseddate\relax\let\accepteddate\relax\let\theasciititle\relax
\let\theasciiauthors\relax
\let\theasciiabstract\relax

\let\theasciiemail\relax
\let\theasciiurl\relax

\ifplaintex
\font\logobig=cmssbx10 scaled 3836
\font\logomed=cmssbx10 scaled 2557
\else
\font\logobig=cmssbx10 scaled 4200
\font\logomed=cmssbx10 scaled 2800
\fi

\long\def\makeagttitle{   
\count0=\startpage
\agt\hfill      
\hbox to 45truept{\vbox to 0pt{\vglue -13truept{\logomed A\kern -.37em{\logobig 
T}\kern -.38em G}\vss}\hss}
\break
{\small Volume \thevolumenumber\ (\thevolumeyear)
\startpage--\finishpage\nl
Published: \publishdate}

\vglue .25truein

{\parskip=0pt\leftskip 0pt plus
1fil\def\\{\par\smallskip}{\Large\bf\thetitle}\par\medskip} \vglue
0.05truein

{\parskip=0pt\leftskip 0pt plus 1fil\def\\{\par}{\sc\theauthors}
\par\medskip}%
 
\vglue 0.03truein 

{\small\leftskip 25truept\rightskip 25truept{\bf Abstract}\stdspace\theabstract

{\bf AMS Classification}\stdspace\theprimaryclass
\ifx\thesecondaryclass\relax\else; \thesecondaryclass\fi\par
{\bf Keywords}\stdspace \thekeywords\par}\vglue 7truept

}   

\ifplaintex
\font\phead=cmsl9 scaled 950
\font\pnum=cmbx10 scaled 913
\font\pfoot=cmsl9 scaled 950
\headline{\vbox to 0pt{\vskip -4.5mm\line{\small\phead\ifnum
\count0=\startpage ISSN 1472-2739 (on-line) 1472-2747 (printed)
\hfill {\pnum\folio}\else\ifodd\count0\def\\{ }%
\ifx\theshorttitle\relax\thetitle\else\theshorttitle\fi\hfill{\pnum\folio}
\else\def\\{ and }{\pnum\folio}\hfill\ifx\theshortauthors\relax\theauthors
\else\theshortauthors\fi\fi\fi}\vss}}
\footline{\vbox to 0pt{\vglue 0mm\line{\small\pfoot\ifnum\count0=\startpage
\copyright\ \gtp\hfill\else
\agt, Volume \thevolumenumber\ (\thevolumeyear)\hfill\fi}\vss}}
\else
\font=cmsl9 scaled 1050
\font\lnum=cmbx10 
\font=cmsl9 scaled 1050
\makeatletter
\def\@oddhead{{\small\lhead\ifnum\count0=\startpage ISSN 1472-2739 
(on-line) 1472-2747 (printed)\hfill {\lnum\number\count0}\else\ifodd\count0
\def\\{ }\ifx\theshorttitle\relax \thetitle \else\theshorttitle\fi\hfill
{\lnum\number\count0}\else\def\\{ and }{\lnum\number\count0}
\hfill\ifx\theshortauthors\relax 
\theauthors\else\theshortauthors\fi\fi\fi}}\def\@evenhead{\@oddhead}
\def\@oddfoot{\small\lfoot\ifnum\count0=\startpage\copyright\ \gtp\hfill\else
\agt, Volume \thevolumenumber\ (\thevolumeyear)\hfill\fi}
\def\@evenfoot{\@oddfoot}
\makeatother
\fi
\let\maketitlepage\makeagttitle
\let\makeshorttitle\maketitlepage
\let\maketitle\maketitlepage

\newwrite\gtoutfile
\long\gdef\makeheadfile{  
{\def\\{, }\def\s{ }
\immediate\openout\gtoutfile head.xxx
\immediate\write\gtoutfile{Proxy-for: \ifx\theasciiauthors\relax
\theauthors\else\theasciiauthors\fi\s<\ifx\theasciiemail\relax\theemail\else\theasciiemail\fi>}
\immediate\write\gtoutfile{\noexpand\\}
\immediate\write\gtoutfile{Authors: \ifx\theasciiauthors\relax
\theauthors\else\theasciiauthors\fi}
{\def\\{ }\immediate\write\gtoutfile{Title: \ifx\theasciititle\relax
\thetitle\else\theasciititle\fi}}
\immediate\write\gtoutfile{Subj-class: GT or SG, GR etc}
\immediate\write\gtoutfile{MSC-class: \theprimaryclass\ifx\thesecondaryclass\relax\else, \thesecondaryclass\fi}
\immediate\write\gtoutfile{Journal-ref: Algebr. Geom. Topol. \thevolumenumber\s
(\thevolumeyear) \startpage-\finishpage}
\immediate\write\gtoutfile{Comments: Published by Algebraic and
Geometric Topology at}
\immediate\write\gtoutfile{\s\s\s  http://www.maths.warwick.ac.uk/agt/AGTVol\thevolumenumber/agt-\thevolumenumber-\thepapernumber.abs.html}
\immediate\write\gtoutfile{\noexpand\\}
\immediate\write\gtoutfile{}
\ifx\theasciiabstract\relax
\immediate\write\gtoutfile{\theabstract}\else
\immediate\write\gtoutfile{\theasciiabstract}\fi
\immediate\write\gtoutfile{}
\immediate\write\gtoutfile{\noexpand\\}
\immediate\write\gtoutfile{}
\immediate\closeout\gtoutfile}}  

\def\maketitlepage{\makeagttitle\makeheadfile}
\let\makeshorttitle\maketitlepage
\let\maketitle\maketitlepage